\documentclass[a4paper,reqno]{amsart}

\usepackage[utf8]{inputenc}
\usepackage{amsmath}
\usepackage{amsfonts}
\usepackage{amssymb}
\usepackage{mathtools}
\usepackage[hidelinks]{hyperref}
\providecommand{\emailauthor}[2]{}
\providecommand{\Newlabel}[2]{}
\numberwithin{equation}{section}
\newtheorem{theorem}{Theorem}[section]
\newtheorem{lemma}[theorem]{Lemma}

\theoremstyle{remark}
\newtheorem{remark}[theorem]{Remark}
\newcommand{\R}{\mathbb R}
\newcommand{\Om}{\Omega}
\newcommand{\norm}[1]{\lVert#1\rVert}
\newcommand{\dd}{\,\mathrm d}

\makeatletter
\@ifundefined{subjclassname@2020}{%
  \@namedef{subjclassname@2020}{\textup{2020} Mathematics Subject Classification}%
}{}
\makeatother

\begin{document}

\title[Boundedness and blow-up in a chemotaxis--haptotaxis system]
{A critical-mass phenomenon in a 2D chemotaxis-haptotaxis system}

\author{Bo Mao}
\address{Institute for Mathematical Sciences, Renmin University of China, Beijing, 100872, China}
\email{maobo@ruc.edu.cn}

\author{Tian Xiang$^*$}
\address{Institute for Mathematical Sciences, Renmin University of China, Beijing, 100872, China}
\email{txiang@ruc.edu.cn}
\thanks{$^*$ Corresponding author.}


\begin{abstract}
We study the Neumann initial-boundary value problem for a minimal
chemotaxis--haptotaxis system with logistic tissue remodeling
\[
\begin{cases}
u_t=\Delta u-\chi\nabla\cdot(u\nabla v)-\xi\nabla\cdot(u\nabla w),\\
v_t=\Delta v-v+u,\\
w_t=-vw+\eta w(1-w)
\end{cases}
\]
in a disk $\Om=B_R(0)\subset\R^2$, with $\chi,\xi,\eta\ge0$
and nonnegative radial initial data $(u_0,v_0,w_0)$.
Under the subcritical mass condition $\chi\|u_0\|_{L^1(\Om)}<8\pi$,
we show that the problem admits a unique global classical solution for every $\eta\ge0$.
If additionally  $w_0\not\equiv0$ and $\eta<\eta^0$ with $\eta^0$ explicitly  depending on $\|u_0\|_{L^1(\Om)}$, $R$, and the support of
$w_0$, we show this solution is uniformly bounded.
Under the same condition, we illustrate that $u$ and $v$ converge exponentially to the
spatial mean $\overline{u_0}=\|u_0\|_{L^1(\Omega)}/|\Omega|$, while $w$ decays
in $W^{1,\infty}(\Om)$ at the optimal rate $\overline{u_0}-\eta$.
The boundedness argument combines radial entropy estimates, auxiliary
elliptic potentials, and a transformation of the $w$-equation.
Fix $\chi>0$ and $\eta\ge0$. For every prescribed mass $m>8\pi/\chi$
and any $B,T_0>0$, we construct smooth positive  initial data
$(u_0,v_0)$ with $\|u_0\|_{L^1(\Omega)}=m$ such that the emanating  solution blows up
before $T_0$ for every admissible radial $w_0$ satisfying $w_0(0)=0$
and $\|\nabla w_0\|_{L^\infty(\Om)}\le B$.
A radial energy estimate and concentration-independent control of $\nabla w$ yield a negative-energy blow-up criterion and the lifespan bound $T_{m,\varepsilon}\le C[\ln(1/\varepsilon)]^{-1/3}$ for all sufficiently small initial concentration scales $\varepsilon$.
Our results extend the $8\pi/\chi$ critical-mass phenomenon of the
chemotaxis-only model ($w\equiv0$) to the chemotaxis-haptotaxis system with tissue remodeling, and give an exponential rate of
convergence in the subcritical regime under the stated condition on $\eta$. Moreover, our results also provide a clearer picture on critical mass phenomena for  the chemotaxis-haptotaxis model in radial setting compared to an earlier  study in planar bounded smooth domains, where solutions are shown to be global when $\chi\|u_0\|_{L^1}<4\pi$,  are bounded when $\eta<\eta_0$ with some $\eta_0<\eta^0$,  and only finite- or infinite-time blow-up alternatives are given when $\chi\|u_0\|_{L^1(\Om)}>4\pi$ and $\eta<\eta_0$.
\end{abstract}

\maketitle

\section{Introduction and statement of main results}

Chemotaxis and haptotaxis describe cell migration in response to
gradients of soluble chemicals and immobile tissue, respectively.
In cancer invasion, cells produce matrix-degrading enzymes, while
degradation and remodeling change the tissue through which they
migrate \cite{ChaplainLolas2005,ChaplainLolas}. The aggregation
mechanism of Keller--Segel systems \cite{KellerSegel,HillenPainter}
then acts in an environment whose spatial structure evolves with
the cells. To deeply examine this interaction without cell proliferation,
we investigate dynamical properties of solutions to a  Neumann initial-boundary problem (IBVP) for the following chemotaxis-haptotaxis model with tissue
remodeling on the disk $\Om=B_R(0)\subset\R^2$, where $R>0$:
\begin{equation}\label{eq:model}
\begin{cases}
u_t=\Delta u-\chi\nabla\cdot(u\nabla v)-\xi\nabla\cdot(u\nabla w),
       &x\in\Om,\ t>0,\\
v_t=\Delta v-v+u, &x\in\Om,\ t>0,\\
w_t=-vw+\eta w(1-w), &x\in\Om,\ t>0,\\
\partial_\nu u=\partial_\nu v=\partial_\nu w=0,
       &x\in\partial\Om,\ t>0,\\
(u,v,w)=(u_0,v_0,w_0), &x\in\Om, \ t=0.
\end{cases}
\end{equation}
Here $u,v,w$ denote the cell density, the signal concentration and the
tissue density, respectively, and $\partial_\nu$ means differentiation with respect to the unit outward  normal $\nu$.
The sensitivities satisfy $\chi,\xi\geq 0$, and $\eta\ge0$ is the remodeling
rate.  The cell mass is conserved
and serves as the principal parameter in the results below.

One  of our starting points is the corresponding Keller--Segel system: 
\begin{equation}\label{eq:KS}
\begin{cases}
u_t=\Delta u-\chi\nabla\cdot(u\nabla v),
       &x\in\Om,\ t>0,\\
v_t=\Delta v-v+u, &x\in\Om,\ t>0,\\
\partial_\nu u=\partial_\nu v=0,
       &x\in\partial\Om,\ t>0,\\
(u,v)=(u_0,v_0), &x\in\Om, \ t=0. 
\end{cases}
\end{equation}
For this  parabolic--parabolic problem, Nagai~et al. \cite{NagaiSenbaYoshida} proved global boundedness below
$4\pi/\chi$ on smooth planar domains and below $8\pi/\chi$ for radial
data on a disk. Entropy estimates play a central role in these
mass-dependent results \cite{GajewskiZacharias}; broader accounts
are given in \cite{HorstmannReview,BellomoReview}. Above the
boundedness threshold, the existence of unbounded solutions and
the finiteness of their lifespan are distinct questions.
Finite-time blow-up constructions for elliptic
\cite{Nagai1995,Nagai2001} and parabolic signals
\cite{HerreroVelazquez}, as well as results allowing finite- or
infinite-time blow-up \cite{HorstmannWang,SenbaSuzuki}, illustrate
this distinction. There are few studies on convergence of bounded solutions, see \cite{Cao-15DCDS, Winkler2010} for small-data convergence to the spatial mean  $(\bar u_0, \bar u_0)$,  convergence to a single equilibrium  \cite{FeireislLaurencotPet}, and convergence to the unique spatial mean under subcritical mass \cite{WangWangYang}. It is natural to ask whether or not mass-critical phenomenon and subcritical mass convergence (and further convergence rate) still  hold in the context of \eqref{eq:model}.

For the same cell--signal--tissue coupling, the closest results are
due to Jin and Xiang~\cite[Theorem 1.1]{JinXiang}, who treated both
elliptic and parabolic signals. On bounded smooth planar domains including smooth convex
domains, they established global classical solvability for every
$\eta\ge0$ under $\chi\norm{u_0}_{L^1(\Om)}<4\pi$.
Uniform boundedness requires a further restriction   on $\eta$ in
terms of the cell mass  and domain diameter. Under this restriction,
their supercritical result, for
$\chi\norm{u_0}_{L^1(\Om)}>4\pi$ outside $4\pi\mathbb N$,
allows either an individual solution unbounded in finite or
infinite time, or approximate blow-up along a family of globally
bounded solutions.
They also proved exponential tissue decay and comparison with
chemotaxis-only dynamics; convergence to a chemotaxis-only solution
in the usual sense requires additional smallness of $\chi$. Comparing these with  the knowledge about the minimal chemotaxis-only system,  we are led to ask whether or not there exists   a rigorous  critical-mass phenomenon for the chemotaxis-haptotaxis model \eqref{eq:model}?

The role of remodeling depends on the accompanying cell kinetics.
Logistic damping in the cell equation can prevent aggregation
\cite{TaoWang2008,TaoWang2009}. Studies combining cell proliferation
with the tissue reaction $\eta w(1-w-u)$ include
\cite{TaoWinkler2014,PangWang2017,KeZheng2018}.
Further boundedness results for models with tissue remodeling
were obtained in \cite{PangWang2018,Jin2018}.
Without remodeling, the mass-dependent boundedness criterion in
\cite{XiangZheng2019} also covers a zero cell source.
In \eqref{eq:model}, cell mass is conserved and remodeling acts
only through $\eta w(1-w)$, so control of aggregation must account
for the evolving tissue drift without cell damping.

Finite-time blow-up in chemotaxis--haptotaxis systems has also been
studied beyond the planar setting. Xiang~\cite{Xiang2022} treated
a higher-dimensional parabolic--elliptic--ODE system, while
Rani and Tyagi~\cite{RaniTyagi2025} attempted its
parabolic--parabolic--ODE counterpart and quasilinear version \cite{RaniTyagi2024} (their blow-up arguments adapted from chemotaxis-only techniques have gaps). These works motivate us to take  a closer
examination of the two-dimensional problem, where logarithmic
concentration and the critical mass govern the energy balance.

The aim of this work is to extend the  $8\pi/\chi$ critical-mass phenomenon and sub-critical mass convergence for the minimal chemotaxis-only system to the corresponding chemotaxis-haptotaxis system  interacting with tissue remodeling. Below this threshold,
we establish global classical solvability for every $\eta\ge0$.
A condition depending on the conserved mass, the disk, and the
region where $w_0$ is positive yields uniform boundedness,
exponential stabilization, and the optimal tissue decay rate.
Above the threshold, for any fixed  remodeling rate $\eta\ge 0$,  we construct finite-time blow-up solutions
at every prescribed mass for  a class of  tissue profiles. The tissue profile is required to
vanish at the center in this construction. Therefore, the radial setting
extends the subcritical mass range and permits a definite lifespan
conclusion for the supercritical data constructed here.

To state our main results more precisely, for fixed $\vartheta\in(0,1)$, we assume
\begin{equation}\label{eq:regular-data}
\begin{gathered}
u_0,v_0,w_0\in C^{2+\vartheta}(\overline\Om),\qquad
u_0,v_0,w_0\ge0,\qquad u_0\not\equiv0,\\
u_0,v_0,w_0\text{ are radial},\qquad
\partial_\nu u_0=\partial_\nu v_0=\partial_\nu w_0=0
\quad\text{on }\partial\Om.
\end{gathered}
\end{equation}
Taking the gradient of the $w$-equation and using $\partial_\nu v|_{\partial\Om}=0$, we obtain
\[
\partial_\nu w(t)|_{\partial\Om}=\partial_\nu w_0|_{\partial\Om}e^{\int_0^t(\eta-v(\tau)-2\eta w(\tau))\dd \tau},
\]
which implies  that $\partial_\nu w_0|_{\partial\Om}=0\Longleftrightarrow \partial_\nu w|_{\partial\Om}=0$. This is why we have conveniently  imposed  homogeneous Neumann boundary data for $w$  (which should not be imposed in principle for a solution of ODE) in \eqref{eq:model}; furthermore, the non-negativity and the Neumann condition  for $w$ imply,
for every $x\in\overline\Om$, that
\[
w_0(x)=0\quad\Longrightarrow\quad
w(x,t)=0,\quad\nabla w(x,t)=0,
\qquad\forall t\in[0,T_m).
\]
Accordingly, pointwise estimates for $w$ and $\nabla w$ below
are proved on $\{w_0>0\}$ and extended by zero without further mention.

We write $\overline{u_0}:=|\Om|^{-1}\int_\Om u_0$ and use
$\operatorname{supp}w_0\subset[0,R]$ for the support of the radial
tissue profile.
Uniform bounds will follow when the accumulated degrading signal
dominates tissue renewal wherever
$w_0$ is positive. A simple sufficient restriction is
$0\le\eta<\norm{u_0}_{L^1(\Om)}/(2\pi R I_1(R))$,
where $I_j$ and $K_j$ denote the respective modified Bessel functions of the
first and second kinds of order $j$, see \eqref{IjKj-def} for their expressions. If $w_0\equiv0$, then $w\equiv0$ and $(u,v)$ solves the
Keller--Segel system \eqref{eq:KS} for every $\eta\ge0$. For $w_0\not\equiv0$, the radial support of $w_0$ gives the
following sufficient condition on $\eta$:
\begin{equation}\label{eq:remodeling-condition}
\eta<\eta^0:=\frac{\norm{u_0}_{L^{1}(\Om)}}{2\pi}
\inf_{r\in(0,R)\cap\operatorname{supp}w_0}
\min\left\{K_0(r)+\frac{K_1(R)}{I_1(R)}I_0(r),\ \
\frac{I_0(r)}{R I_1(R)}\right\}.
\end{equation}
Further, it follows from \eqref{eq:annular-constant} that  $\eta<\overline{u_0}$.

\begin{theorem}
\label{thm:boundedness}
Let $\Om=B_R(0)\subset\R^2$, $\chi,\xi\geq 0$, and let the data satisfy
\eqref{eq:regular-data}. Assume $\chi \norm{u_0}_{L^{1}(\Om)}<8\pi$.
Then, for every fixed $\eta\ge0$, the IBVP  \eqref{eq:model} possesses
a unique global classical solution $(u,v,w)$ with $u,v>0$ and $w\geq 0$ which is locally bounded
in $L^\infty(\Om)\times W^{1,\infty}(\Om)\times W^{1,\infty}(\Om)$.
If, in addition, the smallness 
 \eqref{eq:remodeling-condition} holds, that global  solution is uniformly
bounded in the sense that
\begin{equation}\label{eq:subcritical-uniform}
\sup_{t\ge0}\left(
\norm{u(t)}_{L^{\infty}(\Om)}+\norm{v(t)}_{W^{1,\infty}(\Om)}
+\norm{w(t)}_{W^{1,\infty}(\Om)}\right)<+\infty.
\end{equation}
\end{theorem}

\begin{remark}
Comparing the boundedness for $\chi\|u_0\|_{L^1(\Omega)}<4\pi$ in bounded planar domains \cite[Theorem 1.1 (B1)-(B2)]{JinXiang}, for radial data on a disk, Lemma~\ref{lem:radial-moser} supplies the
additional range $4\pi\le\chi\norm{u_0}_{L^1(\Om)}<8\pi$.
The basic radial remodeling bound also permits a larger $\eta$ than
the disk specialization of \cite[condition (1.7)]{JinXiang}, since
\[
\int_0^\infty\frac{e^{-s-R^2/s}}{4\pi s}\dd s
=\frac{K_0(2R)}{2\pi}<\frac{K_0(R)}{2\pi}<\frac{H(R)}{2\pi}=\frac1{2\pi R I_1(R)};
\]
see the arguments in Lemma ~\ref{lem:radial-resolvent} for the last inequality. Lemma~\ref{lem:radial-resolvent} quantifies the
further improvement when the tissue is supported in an interior annulus. Thus, in the radial setting, we expand the boundedness range  and the range of $\eta$ shown in \cite[Theorem 1.1 (B1)-(B2)]{JinXiang}. 
\end{remark}

Under the same remodeling condition, we demonstrate that the bounded solution
approaches the homogeneous cell--signal equilibrium. Identifying
this limit also determines the precise exponential rate of tissue
decay.

\begin{theorem}\label{thm:decay}
Let $\Om=B_R(0)\subset\R^2$, $\chi,\xi\geq 0$, and let the data satisfy
\eqref{eq:regular-data}.
Assume $\chi \norm{u_0}_{L^{1}(\Om)}<8\pi$
and the smallness \eqref{eq:remodeling-condition}.
Then there exist constants $C,\sigma>0$ such that the global classical solution
obtained in Theorem~\ref{thm:boundedness} satisfies
\[
\norm{u(t)-\overline{u_0}}_{L^{\infty}(\Om)}
+\norm{v(t)-\overline{u_0}}_{W^{1,\infty}(\Om)}
\le Ce^{-\sigma t},\qquad\forall t\ge0, 
\]
and $w$ decays exponentially in an optimal rate: 
\begin{equation}\label{eq:subcritical-tissue-decay}
\norm{w(t)}_{W^{1,\infty}(\Om)}\le C e^{-(\overline{u_0}-\eta)t},
\qquad\forall t\ge0.
\end{equation}
As a result, the system has only the unique equilibrium $(\overline u_0, \overline u_0,0)$ in this setting.

The constants are independent of $t$ but may depend on the data and parameters.
\end{theorem}

\begin{remark}
By the $w$-expression in \eqref{eq:w-exp}, for spatially homogeneous solutions with $u=v=\overline{u_0}$ and
$w_0\equiv b>0$, the tissue component decays exactly at rate
$\overline{u_0}-\eta$ when $\eta<\overline{u_0}$.
Thus, the rate in \eqref{eq:subcritical-tissue-decay} is optimal.

We also would like to add that, in the radial setting, we  not only extend the uniqueness and hence convergence to equilibrium under subcritical mass for the chemotaxis-only system \eqref{eq:KS} \cite[Theorem 1.1]{WangWangYang} to the chemotaxis-haptotaxis system ~\eqref{eq:model}, but also provide an exponential rate for such convergence. 
\end{remark}

For our finite-time blow-up, we consider tissue profiles satisfying
the additional central condition
\begin{equation}\label{eq:tissue}
\begin{gathered}
w_0\in C^{2+\vartheta}(\overline\Om),\qquad
w_0\ge0,\quad w_0\text{ is radial},\\
w_0(0)=0,\qquad \partial_\nu w_0=0\quad\text{on }\partial\Om.
\end{gathered}
\end{equation}
A simple nontrivial example of this class is
$w_0(r)=a[2(r/R)^2-(r/R)^4]$  with $a>0$.
We next show that the supercritical mass blow-up for  \eqref{eq:KS} is locally stable  in the context of  \eqref{eq:model} in some sense. To this end, for $B>0$, we let
\[
\mathcal W_B^0:=
\left\{w_0\text{ satisfying }\eqref{eq:tissue}:
\|\nabla w_0\|_{L^\infty(\Om)}\le B\right\}.
\]

\begin{theorem}
\label{thm:mass}
Let $\Om=B_R(0)\subset\R^2$, $\chi>0$ and $\xi, \eta\ge0$.
For every $B,T_0>0$ and every prescribed
total cell mass $m>8\pi/\chi$, there exist strictly positive
radial functions $u_0,v_0\in C^\infty(\overline\Om)$, chosen independently of $w_0$, satisfying
$\int_\Om u_0=m$, the
Neumann compatibility conditions and
\[
\begin{cases}
-\Delta v_0+v_0=u_0 & \text{in }\Om,\\
\partial_\nu v_0=0 & \text{on }\partial\Om,
\end{cases}
\]
such that, for every $w_0\in\mathcal W_B^0$, the classical solution of
\eqref{eq:model} with initial data $(u_0,v_0,w_0)$ satisfies
\begin{equation}\label{eq:blowup}
0<T_m<T_0,\qquad
\limsup_{t\rightarrow T_m}\norm{u(t)}_{L^\infty(\Om)}=\infty.
\end{equation}
For the family constructed in Section~\ref{sec:blowup}, where
$\varepsilon$ is the spatial concentration scale of $u_0$, there exist $C>0$ independent of $\varepsilon$ and $\varepsilon_0>0$ such that
\begin{equation}\label{eq:lifespan-eps}
T_{m,\varepsilon}\le C[\ln(1/\varepsilon)]^{-1/3}<T_0,
\qquad 0<\varepsilon<\varepsilon_0.
\end{equation}
Here, $T_m$ and $T_{m,\varepsilon}$ are the maximal existence times (lifespan) of \eqref{eq:model} emanating from $(u_0, v_0, w_0)$ and $(u_{0,\varepsilon}, v_{0,\varepsilon}, w_{0})$, respectively.
The family and the constants $C$ and $\varepsilon_0$ can be chosen
independently of $w_0\in\mathcal W_B^0$.
\end{theorem}

\begin{remark}
A precise blow-up criterion and lifespan estimate are detailed in Lemma \ref{prop:criterion}. In radial setting, Theorem \ref{thm:mass} provides a rigorous finite-time blow-up instead of blow-up alternatives provided by \cite[Theorem 1.1, (B3)]{JinXiang} in general bounded domains.   In particular,  we mention that our global existence and  blow-up  hold without any restriction on   $\eta$ and $\xi$, indicating the role of haptotaxis on global existence and blow-up is faint. However,  the blow-up construction uses $w_0(0)=0$, we may not say that the effect of haptotaxis is  entirely negligible.  

When $w\equiv0$, Theorem \ref{thm:mass} recovers the well-known supercritical mass radial blow-up for the minimal Keller-Segel model \eqref{eq:KS} \cite{HerreroVelazquez, HorstmannReview,BellomoReview}. Together with Theorem ~\ref{thm:boundedness}, we not only extend the $8\pi/\chi$ critical-mass phenomenon for \eqref{eq:KS} to \eqref{eq:model}, but also we show in certain sense that the subcritical mass global existence   for  \eqref{eq:KS} is globally stable  in the context of  \eqref{eq:model}, and that the supercritical mass blow-up  for  \eqref{eq:KS} is locally stable  in the context of  \eqref{eq:model}.
\end{remark}
The boundedness theorem imposes no central vanishing on $w_0$,
whereas the blow-up construction uses $w_0(0)=0$.
Uniform boundedness also retains the remodeling condition
\eqref{eq:remodeling-condition}. These results therefore do not give
a complete mass dichotomy for arbitrary tissue data and remodeling
rates.

Our global existence and boundedness rely on the known energy functional
\begin{equation}\label{eq:full-energy}
\mathcal F(u,v,w)=\int_\Om u\ln u-\chi\int_\Om uv
-\xi\int_\Om uw+\frac\chi2\int_\Om(|\nabla v|^2+v^2).
\end{equation}
 Along classical solutions, write
$\mathcal F(t)=\mathcal F(u(t),v(t),w(t))$, the computations in  \cite[Lemma 3.1]{JinXiang} with $\tau=1$ show  that 
\begin{equation}\label{eq:full-energy-identity}
\begin{aligned}
&\frac{d}{dt}\mathcal F(t)
+\int_\Om u|\nabla(\ln u-\chi v-\xi w)|^2
+\chi\int_\Om v_t^2\\
=&\xi\int_\Om uw(v-\eta+\eta w),
\qquad\forall t\in(0,T_m).
\end{aligned}
\end{equation}
The identity extends in integrated form to $t=0$ by continuity
of the initial energy.

For global existence and boundedness, the main difficulty is
to control the tissue gradient without diffusion in the
$w$-equation. We express the accumulated signal as the time
integral of a radial elliptic potential $\Phi$  plus a uniformly bounded
correction $\Psi$ so that $v=\Phi-\Psi_t$ implied by \eqref{eq:elliptic-potentials}. An integrating factor removes the time derivative of
this correction from the tissue equation. Combined with
subcritical entropy control from the subtle radial Trudinger--Moser
inequality (Lemma ~\ref{lem:radial-moser}), the transformed equation yields cell and gradient
bounds on every finite time interval for arbitrary $\eta\ge0$ (Lemmas \ref{lem:subcritical-entropy}-\ref{lem:coupled-bounds}). Semigroup estimates then exclude finite-time blow-up.
Under the restriction  of $\eta$ in  \eqref{eq:remodeling-condition}, a lower bound for the
elliptic potential on the positivity set of $w_0$ gives tissue
decay and damping in the transformed gradient equation.
The coefficient in the entropy inequality becomes integrable,
so the resulting bounds are uniform in time.   In the  derivation of $L^p$ and higher order  gradient bounds, instead of studying the time evolution of $\|u\|_2^2+\|\nabla v\|_4^4+\|\nabla w\|_6^6$ as done in \cite{JinXiang},  we conveniently explore  
\begin{equation*}
\frac{d}{dt}\left(\int_\Om u^2+\int_\Om|\nabla v|^4
+\int_\Om|\nabla(we^{-\Psi})|^6\right).
\end{equation*}
These  make our arguments elegant and greatly  shorten the boundedness proof. 

For stabilization, tissue's $W^{1,\infty}$-decay makes the additional energy
production integrable and permits the identification of stationary
limits (Lemmas \ref{lem:subcritical-compactness}-\ref{lem:subcritical-convergence}). This approach is related to the asymptotic analysis of
invasion models \cite{LitcanuMorales,HillenPainterWinkler} and
convergence to Keller--Segel equilibria
\cite{FeireislLaurencotPet}. The radial uniqueness result of
\cite{WangWangYang} identifies the limit as
$(\overline{u_0},\overline{u_0},0)$. An energy estimate applied to 
\[
\frac d{dt}
\left(\int_\Om(u-\overline{u_0})^2
+\chi\overline{u_0}\int_\Om|\nabla v|^2\right)
\]
yields $(L^2, W^{1,2})$-exponential convergence of $(u-\overline{u_0}, v-\overline{u_0})$. Then we employ Gagliardo--Nirenberg type interpolation inequalities involving  Sobolev norms and H\"older norms (Lemma \ref{GN-Sob-hol}) to lift the convergence to $(L^\infty, W^{1,\infty})$-norm.  Returning to the tissue equation then gives the optimal decay rate
$\overline{u_0}-\eta$ in $W^{1,\infty}$-norm.

For finite-time blow-up, to circumvent the difficulty in digesting  the haptotactic term $\int_\Om uvw$ appearing in \eqref{eq:full-energy-identity} (c.f. ~\cite[Lemma 4.3]{JinXiang} and ~\cite[Lemma 3.6]{RaniTyagi2025}), we use the chemotactic energy
$\mathcal F+\xi\int_\Om uw$. Although this functional is independent
of $w$ (c.f. \eqref{eq:full-energy}), its derivative contains a term involving $\nabla w$ (c.f. \eqref{eq:modified-energy-identity}):
\[
\begin{aligned}
&\frac d{dt}\left(\mathcal F+\xi\int_\Om uw\right)+\chi\int_\Om v_t^2+\int_\Om u|\nabla(\ln u-\chi v)|^2\\
=&\xi\int_\Om u\nabla w\cdot\nabla(\ln u-\chi v), 
\end{aligned}
\]
which must be properly 
controlled to adapt an energy--dissipation argument
\cite{Mizoguchi,Winkler}. Radial integration of the signal equation,
together with $w_0(0)=0$, provides an  $L^\infty$-bound for
$\nabla w$ on each
fixed finite interval that is independent of the concentration
scale (Lemma \ref{lem:tissue}). A crucial  lower bound for the logarithmically weighted entropy of $u$ (Lemma \ref{lem:radial-log-entropy}) is combined with a weighted estimate for $\nabla v$ as well as a delicate rewriting of energy in \eqref{eq:weighted-energy-identity} to  control the negative
part of the chemotactic energy by the $3/4$ power of the
dissipation (Lemma \ref{thm:spatial}): for some $\alpha>0$, 
\[
-\mathcal F(u,v,w)-\xi\int_\Om uw
\le \alpha\left(
1+\chi\int_\Omega v_t^2
+\int_\Om u|\nabla(\ln u-\chi v)|^2
\right)^{\frac34},
\]
leading to a differential inequality with superlinear
growth of order $4/3$ (c.f. \eqref{eq:superlinear-energy}): for some $\beta, \gamma>0$,
\[
\frac d{dt}\left(\mathcal F+\xi\int_\Om uw\right)
\le \beta
-\frac{1}{2\gamma}\left(-\mathcal F-\xi\int_\Om uw\right)_+^{\frac43}.
\]
This enables us to derive a  criterion for finite-time blow-up by adjusting initial energy (Lemma \ref{prop:criterion}).  Finally, concentrating the initial cell
density makes the initial energy tend to $-\infty$ at a logarithmic rate
whenever the prescribed mass exceeds $8\pi/\chi$ (Lemma \ref{lem:initial-energy}). This yields
finite-time blow-up and the lifespan estimate
\eqref{eq:lifespan-eps}.

Section~\ref{sec:boundedness} proves global existence and boundedness,
Section~\ref{sec:decay} establishes stabilization, and
Section~\ref{sec:blowup} constructs finite-time blow-up solutions.

Throughout the paper,  when convenient and no danger of confusion arises, we shall omit integration variables like $\dd x, \dd t$ etc, identify $f=f(t)= f(\cdot,t)$,  and  abbreviate
$\|f\|_p=\|f\|_{L^p(\Om)}$ and
$\|f\|_{s,p}=\|f\|_{W^{s,p}(\Om)}$. Herein  $\Delta_N$ denotes the Neumann Laplacian.
 $\lambda_1$ denotes the first positive
radial Neumann eigenvalue of $-\Delta$ on $B_R$.
 $C>0$ denotes a generic constant. The positive
constants and time-dependent bounds $C_i$ used thereafter are numbered
consecutively in order of first appearance. Unless additional dependence
is displayed, they depend only on the fixed domain, parameters, and
initial data. Time-dependent bounds $C_i(t)$ are chosen non-decreasing
and locally bounded on $[0,\infty)$, independently of the distance
to the maximal existence time. In the concentrating constructions,
independence of the concentration parameter is stated explicitly.

\section{Global existence and uniform boundedness}\label{sec:boundedness}
In this section, we show the global existence and boundedness results
stated in Theorem~\ref{thm:boundedness}. We first recall the local existence
and extensibility criterion from \cite[Lemma 2.1]{JinXiang} (see also \cite{Jin2018}), which are well-established via a proper fixed-point framework and parabolic
regularity theory.
\begin{lemma}\label{lem:local}
Let $\chi,\xi, \eta\ge0$, and let the initial data satisfy
\eqref{eq:regular-data}. There exist a maximal existence time
$T_m\in(0,+\infty]$ and a unique radial classical solution $(u,v,w)$
of~the IBVP \eqref{eq:model} with
\[
(u,v,w)\in
\left[C^0\bigl(\overline\Om\times[0,T_m)\bigr)
\cap C^{2,1}\bigl(\overline\Om\times(0,T_m)\bigr)\right]^3.
\]
The solution satisfies
\[
0<u,\qquad 0<v,\qquad
0\le w\le\max\{1,\norm{w_0}_{L^\infty(\Om)}\}=:K
\quad\text{on }\overline\Om\times(0,T_m).
\]
Moreover, if $T_m<\infty$, then
\begin{equation}\label{eq:local-extensibility}
\limsup_{t\rightarrow T_m}\left(
\norm{u(t)}_{L^\infty(\Om)}+
\norm{v(t)}_{W^{1,\infty}(\Om)}\right)=+\infty.
\end{equation}
\end{lemma}

The following lemma is a consequence of \cite[Lemmas 2.2--2.4]{JinXiang};
its proof is omitted for brevity.
\begin{lemma}\label{lem:mass-signal}
The local-in-time solution satisfies
\[
\|u(t)\|_1=\|u_0\|_1,\qquad\forall t\in[0,T_m).
\]
Moreover, for every $q\in[1,\infty)$ and $r\in[1,2)$,
there exists $C_{1}(q,r)>0$ such that
\[
\|v(t)\|_q+\|v(t)\|_{1,r}\le C_{1}(q,r),
\qquad\forall t\in[0,T_m).
\]
We also have the explicit representation
\begin{equation}\label{eq:w-exp}
w(x,t)=\frac{w_0(x)e^{\eta t-\int_0^t v(x,\tau)\dd \tau}}{1+\eta w_0(x) \int_0^te^{\eta s-\int_0^sv(x,\tau)\dd \tau}\dd s},
\quad\forall(x,t)\in\overline\Om\times[0,T_m),
\end{equation}
and hence
\begin{equation}\label{eq:tissue-direct-bound}
0\le w(x,t)\le w_0(x)e^{\eta t-\int_0^t v(x,s)\dd s},
\quad\forall(x,t)\in\overline\Om\times[0,T_m).
\end{equation}
\end{lemma}

To control $w$ over long times, we seek a lower bound for
the time integral of $v$ in \eqref{eq:tissue-direct-bound}.
For each $t\in[0,T_m)$, introduce the auxiliary elliptic
potentials $\Phi(t)$ and $\Psi(t)$ associated with
$u(t)$ and $v(t)$:
\begin{equation}\label{eq:elliptic-potentials}
\begin{cases}
-\Delta\Phi+\Phi=u,\quad -\Delta\Psi+\Psi=v & \text{in }\Om,\\
\partial_\nu\Phi=\partial_\nu\Psi=0 & \text{on }\partial\Om.
\end{cases}
\end{equation}
We use this notation only in the present section and in
Section~\ref{sec:decay}. The potential $\Phi$ solves the stationary
$v$-equation of~\eqref{eq:model} with $u$ fixed.
Time derivative of the potential $\Psi$ accounts for the difference between $\Phi$
and $v$. Taking $q=4$ in
Lemma~\ref{lem:mass-signal}, elliptic regularity and the maximum
principle yield a constant $C_{2}>0$ such that
\begin{equation}\label{eq:corrector-bounds}
\Psi\ge0,\qquad \|\Psi(t)\|_{1,\infty}\le C_{2},
\qquad\forall t\in[0,T_m).
\end{equation}
Differentiating the equation for $\Psi$ and using  \eqref{eq:elliptic-potentials} and the
$v$-equation of \eqref{eq:model}, we obtain, for each $t\in(0,T_m)$,
\[
\begin{cases}
-\Delta\Psi_t+\Psi_t
=v_t=-\Delta(\Phi-v)+(\Phi-v)
&\text{in }\Om,\\
\partial_\nu\Psi_t=\partial_\nu(\Phi-v)=0
&\text{on }\partial\Om.
\end{cases}
\]
Uniqueness of the Neumann elliptic problem gives
$\Psi_t=\Phi-v$, and hence
\begin{equation}\label{eq:accumulated-signal}
\int_0^t v(x,s)\dd s=\int_0^t\Phi(x,s)\dd s+\Psi(x,0)-\Psi(x,t),\quad\forall (x,t)\in\overline\Om\times[0,T_m).
\end{equation}
Therefore, a lower bound for $\Phi$ controls the time integral of $v$,
up to a uniformly bounded error by \eqref{eq:corrector-bounds}. The following estimate also
records how this lower bound improves on $\operatorname{supp}w_0$.

\begin{lemma}\label{lem:radial-resolvent}
Let $u\ge0$ be smooth and radial on $\overline\Om$, and let $\Phi$ solve
\[
\begin{cases}
-\Delta\Phi+\Phi=u & \text{in }\Om,\\
\partial_\nu\Phi=0 & \text{on }\partial\Om.
\end{cases}
\]
Then, for every $r\in(0,R)$, with $I_j$ and  $K_j$ defined by \eqref{IjKj-def},
\begin{equation}\label{eq:pointwise-resolvent}
\Phi(r)\ge\frac{\|u\|_1}{2\pi}
\min\left\{
K_0(r)+\frac{K_1(R)}{I_1(R)}I_0(r),\ \
\frac{I_0(r)}{R I_1(R)}
\right\}.
\end{equation}
In particular,
\[
\Phi(r)\ge\frac{\|u\|_1}{2\pi R I_1(R)},
\qquad \forall r\in[0,R].
\]
If $w_0\not\equiv0$, then
\begin{equation}\label{eq:annular-constant}
\frac1{R I_1(R)}
\le\inf_{r\in(0,R)\cap\operatorname{supp}w_0}
\min\left\{
K_0(r)+\frac{K_1(R)}{I_1(R)}I_0(r),\
\frac{I_0(r)}{R I_1(R)}
\right\}<\frac{2\pi}{|\Om|}.
\end{equation}
If the radial support of $w_0$ contains $0$ or $R$,
the first inequality in \eqref{eq:annular-constant} is an equality.
\end{lemma}

\begin{proof}
For $r>0$ and $j=0,1$, the standard representations
\cite{Watson}
\begin{equation}\label{IjKj-def}
I_j(r)=\sum_{k=0}^{\infty}\frac{(r/2)^{2k+j}}{k!(k+j)!},
\qquad
K_j(r)=\int_0^\infty e^{-r\cosh t}\cosh(jt)\dd t
\end{equation}
show that $I_0(r)>1$ and $I_1(r),K_0(r),K_1(r)>0$.
To construct the Green representation, we first derive the Wronskian identity.
The derivative formulas $I_0'=I_1$ and $K_0'=-K_1$,
together with the equation
\[
y''+\frac1r y'-y=0,\qquad \forall r>0,
\]
satisfied by both $I_0$ and $K_0$, give
\begin{equation}\label{eq:bessel-wronskian-derivative}
\begin{aligned}
\frac d{dr}\left[r(I_0K_1+I_1K_0)\right]
&=\frac d{dr}\left[r(I_0'K_0-I_0K_0')\right]\\
&=I_0'K_0-I_0K_0'
+r(I_0''K_0-I_0K_0'')\\
&=0,\qquad\forall r>0.
\end{aligned}
\end{equation}
As $r\to0$, the standard asymptotic formulas \cite{Watson}
\[
I_0(r)\sim1,\qquad I_1(r)\sim\frac r2,\qquad
K_0(r)\sim-\ln r,\qquad K_1(r)\sim\frac1r
\]
imply
\[
r\bigl(I_0(r)K_1(r)+I_1(r)K_0(r)\bigr)\to1.
\]
Combining this limit with \eqref{eq:bessel-wronskian-derivative},
we obtain
\[
I_0(r)K_1(r)+I_1(r)K_0(r)=\frac1r,\qquad\forall r>0.
\]
Now set
\[
H(r)=K_0(r)+\frac{K_1(R)}{I_1(R)}I_0(r).
\]
The derivative formulas and the Wronskian identity give
\[
\begin{cases}
 H''+\frac1r H'-H=0,\qquad r>0,\\
 H'(R)=0,\ \
H(R)=\frac1{R I_1(R)}.
\end{cases}
\]
From these  and the derivative formulas for $I_0$ and $K_0$, we obtain the Green representation for the unique solution to the Neumann problem for $\Phi$:
\begin{equation}\label{eq:radial-green}
\Phi(r)=H(r)\int_0^r
I_0(\rho)
u(\rho)\rho\dd\rho+I_0(r)\int_r^RH(\rho)
u(\rho)\rho\dd\rho.
\end{equation}
Observe that
$(rH')'=rH>0$ for $r\in (0,R)$ and $H'(R)=0$, so
$rH'(r)<RH'(R)=0$.
Thus, $H$ is strictly decreasing on $(0,R)$,
while $I_0$ is strictly increasing because $I_0'=I_1>0$. Recalling  $I_0(0)=1$  and $H(R)=1/(R I_1(R))$,
we then deduce from \eqref{eq:radial-green} that
\begin{equation}\label{Phi-est}
 \begin{split}
\Phi(r)&\ge H(r)\int_0^ru(\rho)\rho \dd \rho+I_0(r)H(R)\int_r^Ru(\rho)\rho \dd \rho\\
&\ge H(R)\min\left\{\frac{H(r)}{H(R)}, \ I_0(r)\right\}\int_0^Ru(\rho)\rho \dd \rho\\
&\geq H(R)\int_0^Ru(\rho)\rho \dd \rho=\frac{\|u\|_1}{2\pi R I_1(R)},  \qquad \forall r\in[0,R],
\end{split}
\end{equation}
which along with the formula of $H$ yields \eqref{eq:pointwise-resolvent} and its consequence as well as the first inequality in \eqref{eq:annular-constant} ($\Phi(0)$  and $H(0)$ are  understood in the right  limit sense). Note that
\[
\lim_{r_k\to 0+}\min\left\{\frac{H(r_k)}{H(R)}, \ I_0(r_k)\right\}=1=\lim_{r_k\to R-}\min\left\{\frac{H(r_k)}{H(R)}, \ I_0(r_k)\right\}.
\]
Hence, if  the radial support of $w_0$ contains $0$ or $R$, then
\[
\inf_{r\in(0,R)\cap\operatorname{supp}w_0}
\min\left\{H(r), \ \ I_0(r)H(R)\right\}=H(R)=\frac1{R I_1(R)}.
\]
Finally,
we apply \eqref{eq:radial-green} to $u\equiv1$ to see that $\Phi\equiv1$.
For fixed $r\in(0,R)$, the chain inequalities in \eqref{Phi-est}  show directly that
\[
\frac{2\pi}{|\Omega|}=\frac{2}{R^2}
>\min\left\{H(r), \  \frac{I_0(r)}{R I_1(R)}\right\}.
\]
 This  along with the formula of $H$ entails  the second inequality in  \eqref{eq:annular-constant}.
\end{proof}

We next use the lower bound for $\Phi$ to estimate
the time integral of $v$ and prove exponential decay of $w$.

\begin{lemma}\label{lem:corrector}
Assume $\eta<\eta^0$.
Then there exists $C_3>0$ such that
\begin{equation}\label{eq:corrector-decay}
\|w(t)\|_\infty\le C_3e^{-(\eta^0-\eta)t},
\qquad\forall t\in[0,T_m).
\end{equation}
\end{lemma}

\begin{proof}
On the support of $w_0$, we first use
mass conservation and \eqref{eq:pointwise-resolvent}
to deduce
\[
\Phi\ge\eta^0>\eta.
\]
Combining \eqref{eq:accumulated-signal} with the fact
$\Phi\ge\eta^0$ and $0\le\Psi\le C_2$, we find
\[
\int_0^t v(x,s)\dd s
=\int_0^t\Phi(x,s)\dd s+\Psi(x,0)-\Psi(x,t)
\ge\eta^0t-C_2.
\]
Substituting this lower bound into
\eqref{eq:tissue-direct-bound} yields
\[
0\le w(x,t)\le w_0(x)e^{\eta t-\int_0^t v(x,s)\dd s}\le e^{C_2}w_0(x)e^{-(\eta^0-\eta)t}.
\]
Since $w(x,t)=0$ wherever $w_0(x)=0$,
we extend this bound to all $x\in\overline\Om$.
Taking the supremum over $x\in\overline\Om$,
we obtain \eqref{eq:corrector-decay}.
\end{proof}

To derive the subcritical entropy bound, we first prove
the following exponential inequality of Trudinger–Moser type for radial functions.

\begin{lemma}\label{lem:radial-moser}
For every $\epsilon>0$, there exists $C_4(\epsilon)>0$ such that every radial $v\in H^1(\Omega)$ and every $\alpha>0$ satisfy
\begin{equation}\label{eq:radial-exp}
\ln\int_\Om e^{\alpha v}
\le\left(\frac1{16\pi}+\epsilon\right)\alpha^2
\int_\Om|\nabla v|^2
+\frac{\alpha}{|\Om|}\int_\Om v+C_4(\epsilon).
\end{equation}
\end{lemma}

\begin{proof}
Since $v$ is radial, we have $v-v(R)\in H_0^1(\Omega)$,
where $v(R)$ denotes its constant boundary trace.
If $\|\nabla v\|_2=0$, then $v$ is constant and the estimate
follows immediately. We may therefore assume that $\|\nabla v\|_2>0$.
Then $(v-v(R))/\|\nabla v\|_2$ belongs to $H_0^1(\Omega)$ and
\[
\int_\Om
\left|\nabla\left(\frac{v-v(R)}{\|\nabla v\|_2}\right)\right|^2
=\frac{\int_\Om|\nabla v|^2}{\|\nabla v\|_2^2}
=1.
\]
Hence the Dirichlet Moser inequality~\cite[Theorem 1]{Moser},
with $n=2$, gives a constant $C_5>0$ such that
\begin{equation}\label{eq:radial-moser-quadratic}
\int_\Om
\exp\left(\frac{4\pi(v-v(R))^2}{\|\nabla v\|_2^2}\right)
\le C_5.
\end{equation}
To use \eqref{eq:radial-moser-quadratic}, we first estimate
$\alpha(v-v(R))$ by Young's inequality:
\[
\alpha(v-v(R))
\le\frac{4\pi(v-v(R))^2}{\|\nabla v\|_2^2}
+\frac{\alpha^2\|\nabla v\|_2^2}{16\pi}.
\]
Exponentiating this inequality and integrating over $\Om$,
we deduce from \eqref{eq:radial-moser-quadratic} that
\[
\begin{aligned}
\int_\Om e^{\alpha v}
&=e^{\alpha v(R)}
\int_\Om e^{\alpha(v-v(R))}\\
&\le
e^{\alpha v(R)+\frac{\alpha^2}{16\pi}\|\nabla v\|_2^2}
\int_\Om
\exp\left(\frac{4\pi(v-v(R))^2}{\|\nabla v\|_2^2}\right)\\
&\le
C_5e^{\alpha v(R)+\frac{\alpha^2}{16\pi}\|\nabla v\|_2^2}.
\end{aligned}
\]
Taking logarithms therefore yields
\begin{equation}\label{eq:radial-moser-log}
\ln\int_\Om e^{\alpha v}
\le\alpha v(R)+\frac{\alpha^2}{16\pi}\|\nabla v\|_2^2
+\ln C_5.
\end{equation}
We next compare $v(R)$ with the spatial mean of $v$.
Radial integration and integration by parts give
\[
\begin{aligned}
\frac1{|\Om|}\int_\Om v
&=\frac2{R^2}\int_0^R rv(r)\dd r\\
&=\frac1{R^2}\bigl[r^2v(r)\bigr]_0^R
-\frac1{R^2}\int_0^R r^2v_r(r)\dd r\\
&=v(R)-\frac1{R^2}\int_0^R r^2v_r(r)\dd r.
\end{aligned}
\]
Using this identity, together with Cauchy--Schwarz
and radial integration, we find
\[
\begin{aligned}
\left|v(R)-\frac1{|\Om|}\int_\Om v\right|
&=\frac1{R^2}\left|\int_0^R r^2v_r(r)\dd r\right|\\
&\le\frac1{R^2}
\left(\int_0^R r^3\dd r\right)^{1/2}
\left(\int_0^R r|v_r(r)|^2\dd r\right)^{1/2}\\
&=\frac12\left(\int_0^R r|v_r(r)|^2\dd r\right)^{1/2}\\
&=\frac1{\sqrt{8\pi}}\|\nabla v\|_2.
\end{aligned}
\]
In light of the positivity of $\alpha$, we apply
Young's inequality once more to deduce
\begin{equation}\label{eq:radial-moser-boundary}
\alpha v(R)
\le\frac{\alpha}{|\Om|}\int_\Om v
+\frac{\alpha}{\sqrt{8\pi}}\|\nabla v\|_2\le\frac{\alpha}{|\Om|}\int_\Om v
+\epsilon\alpha^2\|\nabla v\|_2^2
+\frac1{32\pi\epsilon}.
\end{equation}
Finally, substituting \eqref{eq:radial-moser-boundary}
into \eqref{eq:radial-moser-log}, we derive
\eqref{eq:radial-exp}.
\end{proof}

We next use Lemma~\ref{lem:radial-moser} and the energy identity
\eqref{eq:full-energy-identity} to derive the following subcritical
entropy bound for global existence and boundedness. Since standard
entropy estimates cover $\chi=0$, we assume $\chi>0$ in the proof below.

\begin{lemma}\label{lem:subcritical-entropy}
Let $\chi\|u_0\|_1<8\pi$.
There exists $C_6(t)>0$ such that
\begin{equation}\label{eq:entropy-finite}
\int_\Om u\ln(1+u)+\int_\Om v^2+\int_\Om|\nabla v|^2\le C_6(t),\qquad\forall t\in[0,T_m).
\end{equation}
Under the smallness $\eta<\eta^0$ in  \eqref{eq:remodeling-condition}, the bound can be chosen
independently of $t$.
\end{lemma}

\begin{proof}
We first derive a lower bound for the energy $\mathcal F$
defined in \eqref{eq:full-energy}.
Since $\chi\|u_0\|_1<8\pi$, we may choose
$\alpha>\chi$ and $\epsilon>0$ such that
\begin{equation}\label{eq:entropy-parameter-choice}
\frac\chi2-\|u_0\|_1
\left(\frac1{16\pi}+\epsilon\right)\alpha^2>0.
\end{equation}
By mass conservation and Jensen’s inequality, we have
\[
\frac1{\|u_0\|_1}
\int_\Om u\ln\left(\frac{e^{\alpha v}}u\right)
\le
\ln\left(\frac1{\|u_0\|_1}\int_\Om e^{\alpha v}\right).
\]
Rearranging, we obtain
\begin{equation}\label{eq:entropy-jensen}
\int_\Om u\ln u
\ge\alpha\int_\Om uv
+\|u_0\|_1\ln\|u_0\|_1
-\|u_0\|_1\ln\int_\Om e^{\alpha v}.
\end{equation}
Substituting \eqref{eq:entropy-jensen} into the definition
of $\mathcal F$, we apply \eqref{eq:radial-exp}
and use the uniform $L^1$ bound for $v$ from
Lemma~\ref{lem:mass-signal} and the fact $0\le w\le K$
to find a constant $C_7(\epsilon)>0$ such that

\begin{equation}\label{eq:entropy-coercivity}
\begin{aligned}
&\mathcal F=\int_\Om u\ln u-\chi\int_\Om uv-\xi\int_\Om uw
+\frac\chi2\int_\Om v^2+\frac\chi2\int_\Om|\nabla v|^2\\
\ge& (\alpha-\chi)\int_\Om uv+\frac\chi2\int_\Om v^2+\frac\chi2\int_\Om|\nabla v|^2\\
&-\|u_0\|_1\ln\int_\Om e^{\alpha v}+\|u_0\|_1\ln\|u_0\|_1-\xi\|u_0\|_1K\\
\ge&(\alpha-\chi)\int_\Om uv+\frac\chi2\int_\Om v^2+\left[\frac\chi2-\|u_0\|_1\left(\frac1{16\pi}+\epsilon\right)\alpha^2
\right]\int_\Om|\nabla v|^2-C_7(\epsilon).
\end{aligned}
\end{equation}
With $\alpha$ and $\epsilon$ fixed as in
\eqref{eq:entropy-parameter-choice}, we regard $C_7(\epsilon)$
as a fixed constant and denote it simply by $C_7$.
In view of \eqref{eq:entropy-coercivity} and \eqref{eq:entropy-parameter-choice}, we control
$\int_\Om uv$ and
$\int_\Om v^2+\int_\Om|\nabla v|^2$
in terms of $\mathcal F+C_7$.
Returning to the definition of $\mathcal F$ and using
these bounds, together with the fact $0\le w\le K$ and
$s\ln(1+s)-s\ln s\le1$, we deduce, for some $C_8>0$,
\begin{equation}\label{eq:entropy-from-y}
\begin{aligned}
&\int_\Om u\ln(1+u)+\int_\Omega v^2+\int_\Omega|\nabla v|^2\\
\le&\mathcal F+\chi\int_\Om uv+\xi\int_\Om uw+|\Om|
+\left(1-\frac\chi2\right)
\left(\int_\Omega v^2+\int_\Omega|\nabla v|^2\right)\\
\le&C_8(\mathcal F+C_7).
\end{aligned}
\end{equation}
It therefore suffices to bound $\mathcal F+C_7$ from above.
Using the energy identity \eqref{eq:full-energy-identity} and \eqref{eq:entropy-coercivity},
together with the fact that $0\le w\le K$,
we infer, for some $C_9>0$, that
\[
\mathcal F'
\le\xi\|w\|_\infty\int_\Om uv
+\xi\eta\|u_0\|_1\|w\|_\infty^2
\le C_9\|w\|_\infty(\mathcal F+C_7).
\]
We then solve this simple  Gronwall's inequality to deduce
\begin{equation}\label{eq:entropy-gronwall}
\mathcal F(t)+C_7
\le(\mathcal F(0)+C_7)
\exp\left(C_9\int_0^t\|w(s)\|_\infty\dd s\right).
\end{equation}
Substituting $\int_0^t\|w(s)\|_\infty\dd s\le Kt$
into \eqref{eq:entropy-gronwall} and using
\eqref{eq:entropy-from-y}, we conclude that
\[
\int_\Om u\ln(1+u)
+\int_\Om v^2
+\int_\Om|\nabla v|^2
\le C_8(\mathcal F(0)+C_7)e^{C_9Kt},
\]
which proves \eqref{eq:entropy-finite}.
If, in addition, $\eta<\eta^0$, Lemma~\ref{lem:corrector} gives
\[
\int_0^t\|w(s)\|_\infty\dd s
\le C_3\int_0^t e^{-(\eta^0-\eta)s}\dd s
\le\frac{C_3}{\eta^0-\eta}.
\]
The estimate \eqref{eq:entropy-gronwall} therefore provides
a time-independent upper bound for $\mathcal F+C_7$.
We then apply \eqref{eq:entropy-from-y}
to deduce the asserted uniform bound.
\end{proof}

Using the entropy bound from Lemma~\ref{lem:subcritical-entropy},
we next derive a coupled differential estimate for
$\int_\Om u^2+\int_\Om|\nabla v|^4$,
with $\int_\Om|\nabla(we^{-\Psi})|^6$ on the right-hand side.

\begin{lemma}\label{lem:uv-differential}
Let $\chi\|u_0\|_1<8\pi$.
For any $\epsilon>0$, there exists $C_{10}(\epsilon,t)>0$
such that
\begin{equation}\label{eq:uv-coupled-estimate}
\begin{aligned}
&\frac d{dt}\left(\int_\Om u^2+\int_\Om|\nabla v|^4\right)
+\int_\Om|\nabla u|^2
+\frac12\int_\Om|\nabla(|\nabla v|^2)|^2
+4\int_\Om|\nabla v|^4\\
\le&\epsilon\int_\Om|\nabla(we^{-\Psi})|^6
+C_{10}(\epsilon,t),
\qquad\forall t\in(0,T_m).
\end{aligned}
\end{equation}
Under  the smallness $\eta<\eta^0$ in \eqref{eq:remodeling-condition},
$C_{10}(\epsilon,t)$ can be chosen independently of $t$.
\end{lemma}

\begin{proof}
Fix $\epsilon>0$.
Since $\Psi\ge0$, we have
$0\le we^{-\Psi}\le w\le K$. The identity
$\nabla w=e^\Psi\nabla(we^{-\Psi})+w\nabla\Psi$,
together with the uniform $W^{1,\infty}$-bound for $\Psi$
in \eqref{eq:corrector-bounds} and $0\le w\le K$, yields
a constant $C_{11}>0$ such that
\begin{equation}\label{eq:w-gradient-corrector}
|\nabla w|\le C_{11}\bigl(1+|\nabla(we^{-\Psi})|\bigr).
\end{equation}
Testing the $u$-equation with $2u$ and integrating by parts,
we deduce from Young's inequality and
\eqref{eq:w-gradient-corrector} that, for every $\delta>0$,
there exist constants $C_{12}(\delta),C_{13}>0$ such that
\begin{equation}\label{eq:cell-l2-estimate}
\begin{aligned}
\frac d{dt}\int_\Om u^2
=&-2\int_\Om|\nabla u|^2
+2\chi\int_\Om u\nabla u\cdot\nabla v
+2\xi\int_\Om u\nabla u\cdot\nabla w\\
\le&-\frac32\int_\Om|\nabla u|^2
+4\chi^2\int_\Om u^2|\nabla v|^2
+4\xi^2\int_\Om u^2|\nabla w|^2\\
\le&-\frac32\int_\Om|\nabla u|^2
+\delta\int_\Om|\nabla v|^6
+\delta\int_\Om|\nabla(we^{-\Psi})|^6\\
&+C_{12}(\delta)\int_\Om u^3
+C_{13}\int_\Om u^2.
\end{aligned}
\end{equation}
Next, we take the gradient of the $v$-equation and test
the resulting equation with $4|\nabla v|^2\nabla v$.
Using the identity
$2\nabla v\cdot\nabla\Delta v
=\Delta|\nabla v|^2-2|D^2v|^2$
and integrating by parts, upon several computations, we  end up with
\begin{equation}\label{eq:signal-l4-identity}
\begin{aligned}
&\frac d{dt}\int_\Om|\nabla v|^4+2\int_\Om|\nabla(|\nabla v|^2)|^2
+4\int_\Om|\nabla v|^2|D^2v|^2
+4\int_\Om|\nabla v|^4\\
=&-4\int_\Om u\nabla(|\nabla v|^2)\cdot\nabla v
-4\int_\Om u|\nabla v|^2\Delta v
+2\int_{\partial\Om}|\nabla v|^2
\partial_\nu|\nabla v|^2.
\end{aligned}
\end{equation}
The convexity of $\Om$ and the Neumann condition imply
$\partial_\nu|\nabla v|^2\le0$ (indeed, it is zero by radial symmetry of $v$) on $\partial\Om$,
so the boundary term is nonpositive.
Applying $|\Delta v|\le\sqrt2|D^2v|$ and Young's inequality
to the remaining terms, we find, for every $\delta>0$,
a constant $C_{14}(\delta)>0$ such that
\begin{equation}\label{eq:signal-l4-estimate}
\begin{aligned}
&-4\int_\Om u\nabla(|\nabla v|^2)\cdot\nabla v
-4\int_\Om u|\nabla v|^2\Delta v
+2\int_{\partial\Om}|\nabla v|^2
\partial_\nu|\nabla v|^2\\
\le&\int_\Om|\nabla(|\nabla v|^2)|^2
+2\int_\Om|\nabla v|^2|D^2v|^2
+8\int_\Om u^2|\nabla v|^2\\
\le&\int_\Om|\nabla(|\nabla v|^2)|^2
+2\int_\Om|\nabla v|^2|D^2v|^2
+\delta\int_\Om|\nabla v|^6
+C_{14}(\delta)\int_\Om u^3.
\end{aligned}
\end{equation}
Substituting \eqref{eq:signal-l4-estimate} into
\eqref{eq:signal-l4-identity}, we obtain
\begin{equation}\label{eq:signal-l4-differential}
\begin{aligned}
&\frac d{dt}\int_\Om|\nabla v|^4
+\int_\Om|\nabla(|\nabla v|^2)|^2
+2\int_\Om|\nabla v|^2|D^2v|^2
+4\int_\Om|\nabla v|^4\\
\le&\delta\int_\Om|\nabla v|^6
+C_{14}(\delta)\int_\Om u^3.
\end{aligned}
\end{equation}
To estimate the two terms $\int_\Om u^3$ and
$\int_\Om|\nabla v|^6$ on the right-hand sides,
we first use Lemma~\ref{lem:subcritical-entropy} to find
\[
\int_\Om|u\ln u|
\le\int_\Om u\ln(1+u)+\frac{|\Om|}{e}
\le C_6(t)+\frac{|\Om|}{e}.
\]
Using mass conservation and the preceding entropy bound,
we infer from the logarithmic Gagliardo--Nirenberg inequality
\cite{TaoWinkler2014,XiangZheng2019} that, for every $\sigma>0$,
there exists a constant $C_{15}(\sigma,t)>0$ such that
\begin{equation}\label{eq:entropy-interpolation}
\int_\Om u^3
\le\frac{\sigma\|u\ln u\|_1}{C_6(t)+|\Om|/e}
\int_\Om|\nabla u|^2+C_{15}(\sigma,t)
\le\sigma\int_\Om|\nabla u|^2+C_{15}(\sigma,t).
\end{equation}
Similarly, applying the planar Gagliardo--Nirenberg inequality
to $|\nabla v|^2$ and using the fact
$\|\nabla v\|_2^2\le C_6(t)$,
we obtain a constant $C_{16}>0$ such that
\begin{equation}\label{eq:signal-interpolation}
\begin{aligned}
\int_\Om|\nabla v|^6
&=\bigl\||\nabla v|^2\bigr\|_3^3\\
&\le C_{16}\bigl\|\nabla(|\nabla v|^2)\bigr\|_2^2
\bigl\||\nabla v|^2\bigr\|_1
+C_{16}\bigl\||\nabla v|^2\bigr\|_1^3\\
&\le C_{16}C_6(t)\int_\Om|\nabla(|\nabla v|^2)|^2
+C_{16}C_6^3(t).
\end{aligned}
\end{equation}
Adding \eqref{eq:cell-l2-estimate} to
\eqref{eq:signal-l4-differential} together, we estimate the
right-hand side using \eqref{eq:entropy-interpolation}
and \eqref{eq:signal-interpolation}, together with
$\int_\Om u^2\le\int_\Om u^3+|\Om|$.
Thus, for every $\delta,\sigma>0$, there exist constants
$C_{17}(\delta),C_{18}(\delta,t),C_{19}(\delta,\sigma,t)>0$
such that
\[
\begin{aligned}
&\frac d{dt}\left(\int_\Om u^2+\int_\Om|\nabla v|^4\right)
+\frac32\int_\Om|\nabla u|^2
+\int_\Om|\nabla(|\nabla v|^2)|^2\\
&+2\int_\Om|\nabla v|^2|D^2v|^2
+4\int_\Om|\nabla v|^4\\
\le&\delta\int_\Om|\nabla(we^{-\Psi})|^6
+2\delta C_{16}C_6(t)
\int_\Om|\nabla(|\nabla v|^2)|^2
+C_{17}(\delta)\int_\Om u^3+C_{18}(\delta,t)\\
\le&\delta\int_\Om|\nabla(we^{-\Psi})|^6
+2\delta C_{16}C_6(t)
\int_\Om|\nabla(|\nabla v|^2)|^2+\sigma C_{17}(\delta)\int_\Om|\nabla u|^2
+C_{19}(\delta,\sigma,t).
\end{aligned}
\]
We first fix  $\delta>0$ such that $\delta\le\epsilon$ and
$2\delta C_{16}C_6(t)\le1/2$, and then take $\sigma>0$
satisfying $\sigma C_{17}(\delta)\le1/2$.
With these choices, we absorb the two gradient terms into
the left-hand side and discard the nonnegative Hessian term
to derive \eqref{eq:uv-coupled-estimate}.
If, in addition, $\eta<\eta^0$,
Lemma~\ref{lem:subcritical-entropy} provides a time-independent
entropy bound $C_6$.
The parameters $\delta$ and $\sigma$ can then be fixed
independently of $t$, and the remainders
$C_{15}(\sigma,t)$ and $C_{18}(\delta,t)$,
and hence $C_{19}(\delta,\sigma,t)$,
can be chosen uniformly bounded in time.
Hence \eqref{eq:uv-coupled-estimate} holds with
a time-independent constant $C_{10}(\epsilon)$.
\end{proof}

To control the remaining term in \eqref{eq:uv-coupled-estimate},
we next derive a differential inequality for
$\int_\Om|\nabla(we^{-\Psi})|^6$ in terms of $\int_\Om u^3$.

\begin{lemma}
There exists $C_{20}>0$ such that,  for all $t\in(0,T_m)$, 
\begin{equation}\label{eq:h-six-finite}
\frac d{dt}\int_\Om|\nabla(we^{-\Psi})|^6
\le C_{20}\left(
\int_\Om|\nabla(we^{-\Psi})|^6+\int_\Om u^3+1
\right).
\end{equation}
Under  the smallness $\eta<\eta^0$ in \eqref{eq:remodeling-condition}, there exists $C_{21}>0$
such that, for all $t\in(0,T_m)$, 
\begin{equation}\label{eq:h-six-uniform}
\frac d{dt}\int_\Om|\nabla(we^{-\Psi})|^6
+3(\eta^0-\eta)\int_\Om|\nabla(we^{-\Psi})|^6
\le C_{21}\left(\int_\Om u^3+1\right).
\end{equation}
\end{lemma}

\begin{proof}
By the elliptic equations defining $\Phi$ and $\Psi$ in
\eqref{eq:elliptic-potentials} and the $v$-equation
of~\eqref{eq:model}, we have $v=\Phi-\Psi_t$.
Using this identity, the $w$-equation of~\eqref{eq:model}, and
$\nabla w=e^\Psi\nabla(we^{-\Psi})+w\nabla\Psi$,
we compute

\begin{equation}\label{eq:h-gradient-equation}
\begin{aligned}
&\partial_t\nabla(we^{-\Psi})\\
=&\nabla\bigl(e^{-\Psi}(w_t-w\Psi_t)\bigr)\\
=&-\nabla\bigl((v-\eta+\eta w+\Psi_t)we^{-\Psi}\bigr)\\
=&-\nabla\bigl((\Phi-\eta+\eta w)we^{-\Psi}\bigr)\\
=&-(\Phi-\eta+\eta w)\nabla(we^{-\Psi})
-we^{-\Psi}\nabla\Phi-\eta we^{-\Psi}
\bigl[e^\Psi\nabla(we^{-\Psi})+w\nabla\Psi\bigr]\\
=&-(\Phi-\eta+\eta w)\nabla(we^{-\Psi})
-we^{-\Psi}\nabla\Phi-\eta w\nabla(we^{-\Psi})
-\eta w^2e^{-\Psi}\nabla\Psi\\
=&-(\Phi-\eta+2\eta w)\nabla(we^{-\Psi})
-we^{-\Psi}\nabla\Phi-\eta w^2e^{-\Psi}\nabla\Psi.
\end{aligned}
\end{equation}
Next, using elliptic regularity in \eqref{eq:elliptic-potentials}  and interpolation, we find, for some
constant $C_{22}>0$,  that
\begin{equation}\label{eq:z-gradient-six}
\|\nabla\Phi\|_6^6
\le C_{22}\|u\|_{3/2}^6
\le C_{22}\|u\|_1^3\|u\|_3^3
=C_{22}\|u_0\|_1^3\int_\Om u^3.
\end{equation}
Taking the inner product of \eqref{eq:h-gradient-equation}
with $6|\nabla(we^{-\Psi})|^4\nabla(we^{-\Psi})$
and integrating over $\Om$, we use the nonnegativity of $\Psi$
and its uniform $W^{1,\infty}$-bound in
\eqref{eq:corrector-bounds}, together with
\eqref{eq:z-gradient-six}, $0\le w\le K$,
and Young's inequality, to obtain
\begin{equation}\label{eq:h-six-identity}
\begin{aligned}
&\frac d{dt}\int_\Om|\nabla(we^{-\Psi})|^6\\
=&-6\int_\Om(\Phi-\eta+2\eta w)|\nabla(we^{-\Psi})|^6-6\int_\Om we^{-\Psi}|\nabla(we^{-\Psi})|^4
\nabla(we^{-\Psi})\cdot\nabla\Phi\\
&-6\eta\int_\Om w^2e^{-\Psi}|\nabla(we^{-\Psi})|^4
\nabla(we^{-\Psi})\cdot\nabla\Psi\\
\le&6\eta\int_\Om|\nabla(we^{-\Psi})|^6
+6K\int_\Om|\nabla(we^{-\Psi})|^5|\nabla\Phi|+6\eta K^2C_2\int_\Om|\nabla(we^{-\Psi})|^5\\
\le& C_{20}\left(
\int_\Om|\nabla(we^{-\Psi})|^6+\int_\Om u^3+1
\right).
\end{aligned}
\end{equation}
We have thus proved \eqref{eq:h-six-finite}.
Under the condition $\eta<\eta^0$ in  \eqref{eq:remodeling-condition}, we recall from
the proof of Lemma~\ref{lem:corrector} that
$\Phi-\eta\ge\eta^0-\eta>0$ on the radial support of $w_0$.
Since $\nabla(we^{-\Psi})=0$ outside this support
and $w\ge0$, we use the equality in
\eqref{eq:h-six-identity}, together with Young's inequality
and \eqref{eq:z-gradient-six}, to obtain
\[
\begin{aligned}
&\frac d{dt}\int_\Om|\nabla(we^{-\Psi})|^6
+6(\eta^0-\eta)\int_\Om|\nabla(we^{-\Psi})|^6\\
\le&\frac d{dt}\int_\Om|\nabla(we^{-\Psi})|^6
+6\int_\Om(\Phi-\eta+2\eta w)|\nabla(we^{-\Psi})|^6\\
=&-6\int_\Om we^{-\Psi}|\nabla(we^{-\Psi})|^4
\nabla(we^{-\Psi})\cdot\nabla\Phi\\
&-6\eta\int_\Om w^2e^{-\Psi}|\nabla(we^{-\Psi})|^4
\nabla(we^{-\Psi})\cdot\nabla\Psi\\
\le&6K\int_\Om|\nabla(we^{-\Psi})|^5|\nabla\Phi|
+6\eta K^2C_2\int_\Om|\nabla(we^{-\Psi})|^5\\
\le&3(\eta^0-\eta)\int_\Om|\nabla(we^{-\Psi})|^6
+C_{21}\left(\int_\Om u^3+1\right).
\end{aligned}
\]
Finally, absorbing the first term on the right-hand side
into the left-hand side, we arrive at
\eqref{eq:h-six-uniform}.
\end{proof}

 Now, we combine the preceding two lemmas and apply Gronwall's
inequality to bound $\|u\|_2$, $\|\nabla v\|_4$, and
$\|\nabla(we^{-\Psi})\|_6$.
If $\eta<\eta^0$, these bounds are uniform in time.

\begin{lemma}\label{lem:coupled-bounds}
Let $\chi\|u_0\|_{1}<8\pi$. There exists $C_{23}(t)>0$ such that
\begin{equation}\label{eq:coupled-finite}
\|u\|_{2}+\|\nabla v\|_{4}
+\|\nabla(we^{-\Psi})\|_{6}\le C_{23}(t),
\quad\forall t\in[0,T_m).
\end{equation}
Under  the smallness $\eta<\eta^0$ in \eqref{eq:remodeling-condition}, the bound can be chosen
independently of $t$.
\end{lemma}
\begin{proof}
Define  
\[
y(t):=\int_\Om u^2+\int_\Om|\nabla v|^4
+\int_\Om|\nabla(we^{-\Psi})|^6.
\]
We apply Lemma~\ref{lem:uv-differential} with $\epsilon=1$
and add \eqref{eq:uv-coupled-estimate} to
\eqref{eq:h-six-finite}.
Choosing $\sigma>0$ such that $C_{20}\sigma\le1/2$
and using the interpolation estimate
\eqref{eq:entropy-interpolation}, we obtain
a constant $C_{24}(t)>0$ such that
\[
\begin{aligned}
&y'(t)+\int_\Om|\nabla u|^2
+\frac12\int_\Om|\nabla(|\nabla v|^2)|^2
+4\int_\Om|\nabla v|^4\\
\le& (1+C_{20})\int_\Om|\nabla(we^{-\Psi})|^6
+C_{20}\int_\Om u^3+C_{10}(1,t)+C_{20}\\
\le& C_{20}\sigma\int_\Om|\nabla u|^2
+(1+C_{20})\int_\Om|\nabla(we^{-\Psi})|^6+C_{24}(t)\\
\le&\frac12\int_\Om|\nabla u|^2
+(1+C_{20})\int_\Om|\nabla(we^{-\Psi})|^6+C_{24}(t).
\end{aligned}
\]
Absorbing the first term on the right-hand side into
the left-hand side and discarding the remaining nonnegative
terms, we obtain, by the definition of $y$,
\[
y'(t)\le(1+C_{20})y(t)+C_{24}(t).
\]
We then solve this simple Gronwall's inequality to conclude that
\[
y(t)\le y(0)e^{(1+C_{20})t}+\int_0^t e^{(1+C_{20})(t-s)}C_{24}(s)\,ds.
\]
We have thus shown \eqref{eq:coupled-finite}.
Under the condition $\eta<\eta^0$, the remainder constants in
\eqref{eq:uv-coupled-estimate} and
\eqref{eq:entropy-interpolation} can be chosen independently of $t$.
We fix $0<\epsilon<3(\eta^0-\eta)/2$ and write
$C_{10}$ for $C_{10}(\epsilon)$.
Adding \eqref{eq:uv-coupled-estimate} to
\eqref{eq:h-six-uniform} and absorbing
$\epsilon\int_\Om|\nabla(we^{-\Psi})|^6$
into the damping term, we arrive at
\begin{equation}\label{eq:coupled-damped}
\begin{aligned}
&y'(t)+\int_\Om|\nabla u|^2
+\frac12\int_\Om|\nabla(|\nabla v|^2)|^2\\
&+4\int_\Om|\nabla v|^4
+\frac{3(\eta^0-\eta)}2
\int_\Om|\nabla(we^{-\Psi})|^6\le C_{21}\int_\Om u^3+C_{10}+C_{21}.
\end{aligned}
\end{equation}
Fix $\sigma>0$ such that $(C_{21}+1)\sigma\le1/2$.
Using the fact that $u^2\le u^3+1$ and \eqref{eq:entropy-interpolation},
we find a constant $C_{25}>0$ such that
\begin{equation}\label{eq:coupled-absorption}
\begin{aligned}
C_{21}\int_\Om u^3+\int_\Om u^2
&\le(C_{21}+1)\int_\Om u^3+|\Om|\\
&\le(C_{21}+1)\sigma\int_\Om|\nabla u|^2
+(C_{21}+1)C_{15}(\sigma)+|\Om|\\
&\le\frac12\int_\Om|\nabla u|^2+C_{25}.
\end{aligned}
\end{equation}
Adding $\int_\Om u^2$ to both sides of
\eqref{eq:coupled-damped} and using
\eqref{eq:coupled-absorption},
we obtain, for some  $C_{26}>0$,
\[
\begin{aligned}
&y'(t)+\int_\Om u^2
+\frac12\int_\Om|\nabla u|^2
+\frac12\int_\Om|\nabla(|\nabla v|^2)|^2\\
&+4\int_\Om|\nabla v|^4
+\frac{3(\eta^0-\eta)}2
\int_\Om|\nabla(we^{-\Psi})|^6
\le C_{26},
\end{aligned}
\]
which along with  the definition
of $y$ leads to
\[
y'(t)+\min\left\{1,\frac{3(\eta^0-\eta)}2\right\}y(t)
\le C_{26}.
\]
Since $\eta<\eta^0$, we deduce from Gronwall's inequality
that $y$ is uniformly bounded in time.
Hence, the  desired uniform bound  \eqref{eq:coupled-finite} holds.
\end{proof}

Using the bounds from Lemma~\ref{lem:coupled-bounds}
and the well-known  $L^p$-$L^q$ Neumann semigroup estimates \cite{Winkler2010, Cao-15DCDS}, we now derive
$L^\infty(\Om)$ bounds for $u$ and
$W^{1,\infty}(\Om)$ bounds for $v$ and $w$. We offer these standard arguments for completeness.

\begin{lemma}\label{lem:bounded-bootstrap}
Let $\chi\|u_0\|_{1}<8\pi$. If $T_m<\infty$, there exists
$C_{27}(T_m)>0$ such that
\[
\|u(t)\|_\infty+\|v(t)\|_{1,\infty}+\|w(t)\|_{1,\infty}
\le C_{27}(T_m),\qquad\forall t\in[0,T_m).
\]
Under  the smallness $\eta<\eta^0$ in \eqref{eq:remodeling-condition}, these norms are uniformly
bounded on $(0,T_m)$.
\end{lemma}

\begin{proof}
We first assume $T_m<\infty$ and derive estimates valid
for all $t\in[0,T_m)$.
By the triangle inequality, \eqref{eq:w-gradient-corrector},
H\"older's inequality, and Lemma~\ref{lem:coupled-bounds},
we obtain a constant $C_{28}(T_m)>0$ such that
\[
\begin{aligned}
\|\chi\nabla v+\xi\nabla w\|_4
&\le\chi\|\nabla v\|_4+\xi\|\nabla w\|_4\\
&\le\chi\|\nabla v\|_4
+\xi C_{11}\left(
|\Om|^{1/4}+\|\nabla(we^{-\Psi})\|_4
\right)\\
&\le\chi\|\nabla v\|_4
+\xi C_{11}\left(
|\Om|^{1/4}
+|\Om|^{1/12}\|\nabla(we^{-\Psi})\|_6
\right)\\
&\le C_{28}(T_m).
\end{aligned}
\]
Together with the $L^2$ bound for $u$ from
Lemma~\ref{lem:coupled-bounds}, this estimate allows us
to bound $u(\chi\nabla v+\xi\nabla w)$ in $L^{4/3}$.
Applying H\"older's inequality and the Neumann semigroup
estimates for the $u$-equation, we find constants
$C_{29},C_{30}(T_m)>0$ such that
\[
\begin{aligned}
\|u\|_3&\le C_{29}\|u_0\|_3\\
&+C_{29}\int_0^t
\bigl(1+(t-s)^{-11/12}\bigr)e^{-\lambda_1(t-s)}
\|u(s)(\chi\nabla v(s)+\xi\nabla w(s))\|_{4/3}\dd s\\
&\le C_{29}\|u_0\|_3\\
&+C_{29}\int_0^t
\bigl(1+(t-s)^{-11/12}\bigr)e^{-\lambda_1(t-s)}
\|u(s)\|_2\|\chi\nabla v(s)+\xi\nabla w(s)\|_4\dd s\\
\le& C_{29}\|u_0\|_3+C_{29}C_{23}(T_m)C_{28}(T_m)
\int_0^\infty(1+s^{-11/12})e^{-\lambda_1s}\dd s\\
\le&C_{30}(T_m).
\end{aligned}
\]
We use the preceding $L^3$-bound for $u$ and parabolic
regularity for the $v$-equation of \eqref{eq:model} to find
a constant $C_{31}(T_m)>0$ such that
$\|v(t)\|_{1,\infty}\le C_{31}(T_m)$
for all $t\in[0,T_m)$. Rewriting the $w$-equation as   $w_t=-(v-\eta+\eta w)w$, we then infer
\begin{equation}\label{w-gradw-sol}
w(t)=w(s)e^{-\int_s^t (v-\eta+\eta w)(\tau) d\tau},\ \
\  \  \partial_t(\nabla w)=-(v-\eta+2\eta w)\nabla w-w\nabla v.
\end{equation}
These along with $w, \eta\geq 0$ enable us to further deduce
\begin{equation}\label{gradw-bdd}
\begin{split}
|\nabla w(t)|&\le |\nabla w_0|e^{-\int_0^t(v-\eta+2\eta w)(\tau)d\tau}+\int_0^te^{-\int_s^t(v-\eta+2\eta w)(\tau)d\tau}w(s)|\nabla v(s)|ds\\
&\le |\nabla w_0|e^{\eta t-\int_0^t v(s)\dd s}
+w(t)\int_0^t|\nabla v(s)|\dd s.
\end{split}
\end{equation}
Since $v\ge0$ and $0\le w\le K$, we deduce from the preceding
estimates there exists a constant $C_{32}(T_m)>0$
such that
\[
\|w(t)\|_{1,\infty}
\le K+e^{\eta T_m}\|\nabla w_0\|_\infty
+KT_mC_{31}(T_m)\le C_{32}(T_m).
\]
Returning to the $u$-equation of~\eqref{eq:model}, writing $f=\chi u\nabla v+\xi u\nabla w$, with $L^3$-bound of $f$, we obtain
constants $C_{33}, C_{34}(T_m),C_{35}(T_m)>0$ such that
\begin{equation}\label{eq:bootstrap-maximum}
\begin{aligned}
\|u(t)\|_\infty\le&\|u_0\|_\infty+C_{33}\int_0^t
\bigl(1+(t-s)^{-5/6}\bigr)e^{-\lambda_1(t-s)}
\|f(s)\|_3\dd s\\
\le&\|u_0\|_\infty+C_{34}(T_m)\int_0^\infty(1+s^{-5/6})e^{-\lambda_1s}\dd s\\
\le&C_{35}(T_m).
\end{aligned}
\end{equation}
Finally, we combine the bounds for $u$, $v$, and $w$
to finish the first assertion.
Under the smallness condition $\eta<\eta^0$, we employ
Lemma~\ref{lem:coupled-bounds} in the preceding estimates
to obtain uniform bounds for $\|u(t)\|_3$ and
$\|v(t)\|_{1,\infty}$.
We now return to the pointwise estimate for $\nabla w$.
The uniform bound for $\|\nabla v(t)\|_\infty$,
together with \eqref{gradw-bdd} and the exponential decay estimates established in the proof
of Lemma~\ref{lem:corrector}, yields a uniform bound for
$\|w(t)\|_{1,\infty}$, since
$(1+t)e^{-(\eta^0-\eta)t}$ is bounded.
Substituting these bounds into \eqref{eq:bootstrap-maximum},
we conclude that $\|u(t)\|_\infty$ is uniformly bounded.
We have thus proved the second assertion.
\end{proof}

\begin{proof}[Proof of Theorem~\ref{thm:boundedness}]
We consider the unique maximal classical solution $(u,v,w)$
provided by Lemma~\ref{lem:local}.
If $T_m<\infty$, we deduce from
Lemma~\ref{lem:bounded-bootstrap} that the norms in
\eqref{eq:local-extensibility} remain bounded as $t\rightarrow T_m$.
We thus contradict the extensibility criterion
and conclude that the solution exists globally.
By the regularity of the global solution, we also obtain
the asserted local boundedness.
If $\eta<\eta^0$, we deduce \eqref{eq:subcritical-uniform}
directly from the uniform estimate in
Lemma~\ref{lem:bounded-bootstrap}.
\end{proof}

\section{Exponential convergence to equilibrium}\label{sec:decay}
Throughout this section, we assume that $\chi\|u_0\|_1<8\pi$ and
$\eta<\eta^0$, so that
Theorem~\ref{thm:boundedness} yields a global uniformly bounded
solution.
We prove Theorem~\ref{thm:decay} using energy dissipation,
compactness, and exponential decay estimates.
The case $\chi=0$ follows from standard energy and semigroup
estimates (See \cite[Theorem 1.1]{Xiang2022} for the case of  $\eta=0$). We therefore give the detailed proof only for $\chi>0$.

\begin{lemma}\label{lem:subcritical-compactness}
There exists $C_{36}>0$ such that
\begin{equation}\label{eq:preliminary-tissue-decay}
\|w(t)\|_{1,\infty}\le C_{36}(1+t)e^{-(\eta^0-\eta)t},
\qquad\forall t\ge0.
\end{equation}
Moreover, for every sequence $t_j\to\infty$,
we can extract a subsequence, still denoted by $t_j$,
and find nonnegative radial functions
$(u_\infty,v_\infty)\in C(\overline\Om)\times C^1(\overline\Om)$
such that
\begin{equation}\label{eq:compact-spatial}
\sup_{t_j\le t\le t_j+1}
\left(
\|u(t)-u_\infty\|_\infty
+\|v(t)-v_\infty\|_{1,\infty}
\right)\to0
\qquad\text{as }j\to\infty.
\end{equation}
\end{lemma}

\begin{proof}
By the exponential decay estimate for $w$ in
\eqref{eq:corrector-decay}, we have
$\|w\|_\infty\le C_3e^{-(\eta^0-\eta)t}$.
Utilizing the lower bound
$\int_0^t v(x,s)\dd s\ge\eta^0t-C_2$
on $\operatorname{supp}w_0$, as established in the proof of
Lemma~\ref{lem:corrector}, we  infer from \eqref{gradw-bdd} that
\[
\begin{aligned}
\|\nabla w\|_\infty
&\le
\left\|e^{\eta t-\int_0^t v(\cdot,s)\dd s}
\nabla w_0\right\|_\infty
+\|w\|_\infty\int_0^t\|\nabla v(s)\|_\infty\dd s\\
&\le e^{C_2}\|\nabla w_0\|_\infty e^{-(\eta^0-\eta)t}
+C_3e^{-(\eta^0-\eta)t}
\int_0^t\|\nabla v(s)\|_\infty\dd s\\
&\le\left(
e^{C_2}\|\nabla w_0\|_\infty
+C_3t\sup_{s\ge0}\|\nabla v(s)\|_\infty
\right)e^{-(\eta^0-\eta)t}.
\end{aligned}
\]
Since $\nabla v$ is uniformly bounded, we combine this
gradient estimate with the exponential decay estimate for $w$
in \eqref{eq:corrector-decay} to conclude
\eqref{eq:preliminary-tissue-decay}.

Next, using the exponential  decay of $w$ and the uniform bounds for $v$ and $w$,
we can choose a constant $C_{37}>0$ such that
\[
\left|\xi\int_\Om uw(v-\eta+\eta w)\right|\le\xi\|u_0\|_1
\bigl(\|v\|_\infty+\eta+\eta\|w\|_\infty\bigr)
\|w\|_\infty\le C_{37}e^{-(\eta^0-\eta)t}.
\]
For any $T>1$, we integrate the energy identity
\eqref{eq:full-energy-identity} over $(1,T)$ and use
the preceding estimate, together with
$\mathcal F(T)\ge-C_7$ from \eqref{eq:entropy-coercivity},
to obtain
\[
\begin{aligned}
&\int_1^T\int_\Om
\left(u|\nabla(\ln u-\chi v-\xi w)|^2+\chi v_t^2\right)\dd x\dd t\\
=&\mathcal F(1)-\mathcal F(T)
+\xi\int_1^T\int_\Om uw(v-\eta+\eta w)\dd x\dd t\\
\le&\mathcal F(1)+C_7
+C_{37}\int_1^T e^{-(\eta^0-\eta)t}\dd t\\
\le&\mathcal F(1)+C_7+\frac{C_{37}}{\eta^0-\eta}.
\end{aligned}
\]
Since the right-hand side is independent of $T$,
we let $T\to\infty$ and conclude that
\begin{equation}\label{eq:finite-total-dissipation}
\int_1^\infty\int_\Om
\left(u|\nabla(\ln u-\chi v-\xi w)|^2+\chi v_t^2\right)\dd x\dd t<\infty.
\end{equation}
Moreover, we apply the standard parabolic H\"older estimates to the $u$-equation
and parabolic $L^8$ estimates to the $v$-equation
on each interval $[t-1,t+1]$.
Combining these estimates with the uniform bounds and
the parabolic Sobolev embedding
\cite[Chapters IV, VI and VII]{Lieberman},
we find $\alpha\in(0,1/2)$ and a constant $C_{38}>0$ such that
\begin{equation}\label{eq:uniform-holder-bounds}
\|u(t)\|_{C^\alpha(\overline\Om)}
+\|v(t)\|_{C^{1+\alpha}(\overline\Om)}
\le C_{38},\qquad\forall t\ge1.
\end{equation}
Given any sequence $t_j\to\infty$, we apply
the Arzel\`a--Ascoli theorem and
\eqref{eq:uniform-holder-bounds} to extract a subsequence,
still denoted by $t_j$, and find nonnegative radial functions
$(u_\infty,v_\infty)\in C(\overline\Om)\times C^1(\overline\Om)$
such that
\[
\|u(t_j)-u_\infty\|_\infty
+\|v(t_j)-v_\infty\|_{1,\infty}\to0.
\]
We then show that this convergence is uniform on $[t_j,t_j+1]$
by controlling the time variation of $u$ and $v$.
By \eqref{eq:finite-total-dissipation}, we have, as $j\to \infty$,
\begin{equation}\label{eq:vanishing-tail-dissipation}
\int_{t_j}^{t_j+1}\int_\Om
\left(u|\nabla(\ln u-\chi v-\xi w)|^2+\chi v_t^2\right)\dd x\dd \tau
\to 0.
\end{equation}
By the mass conservation  $\|u\|_1=\|u_0\|_1$ and   the Cauchy--Schwarz inequality, we estimate the flux as follows:
\begin{equation}\label{eq:vanishing-flux}
\begin{aligned}
&\int_{t_j}^{t_j+1}
\|\nabla u-\chi u\nabla v-\xi u\nabla w\|_1\dd\tau\\
\le&\|u_0\|_1^{1/2}
\left(\int_{t_j}^{t_j+1}\int_\Om
u|\nabla(\ln u-\chi v-\xi w)|^2\dd x\dd\tau\right)^{1/2}\\
\le&\|u_0\|_1^{1/2}
\left(\int_{t_j}^{t_j+1}\int_\Om
\left(u|\nabla(\ln u-\chi v-\xi w)|^2+\chi v_t^2\right)
\dd x\dd\tau\right)^{1/2}\to0.
\end{aligned}
\end{equation}
For any $\varphi\in C^1(\overline\Om)$, we test
the $u$-equation of~\eqref{eq:model} against $\varphi$
and integrate by parts to infer
\[
\frac d{d\tau}\int_\Om u(\tau)\varphi
=-\int_\Om
(\nabla u-\chi u\nabla v-\xi u\nabla w)\cdot\nabla\varphi.
\]
Integrating from $t_j$ to $t\in[t_j,t_j+1]$
and applying H\"older's inequality, we obtain
\[
\begin{aligned}
\left|\int_\Om(u(t)-u(t_j))\varphi\right|
&=\left|\int_{t_j}^{t}\int_\Om
(\nabla u-\chi u\nabla v-\xi u\nabla w)
\cdot\nabla\varphi\dd x\dd\tau\right|\\
&\le\|\nabla\varphi\|_\infty
\int_{t_j}^{t_j+1}
\|\nabla u-\chi u\nabla v-\xi u\nabla w\|_1\dd\tau.
\end{aligned}
\]
Taking the supremum over $t\in[t_j,t_j+1]$
and using \eqref{eq:vanishing-flux}, we conclude that
\[
\sup_{t_j\le t\le t_j+1}
\left|\int_\Om(u(t)-u(t_j))\varphi\right|\to0.
\]
For $v$, we similarly integrate in time and apply
the Cauchy--Schwarz inequality to write, for every
$t\in[t_j,t_j+1]$,
\[
\|v(t)-v(t_j)\|_2=\left\|\int_{t_j}^{t}v_t(\tau)\dd\tau\right\|_2\le\int_{t_j}^{t}\|v_t(\tau)\|_2\dd\tau\le\left(\int_{t_j}^{t_j+1}
\|v_t(\tau)\|_2^2\dd\tau\right)^{1/2}.
\]
Since $\chi>0$, we know from
\eqref{eq:vanishing-tail-dissipation} that the final bound
tends to zero.
Taking the supremum over $[t_j,t_j+1]$, we therefore see
\[
\sup_{t_j\le t\le t_j+1}\|v(t)-v(t_j)\|_2\to0.
\]
We finally combine these time-variation estimates with
spatial compactness.
For any $s_j\in[t_j,t_j+1]$, we use
\eqref{eq:uniform-holder-bounds} to extract a subsequence of
$(u(s_j),v(s_j))$ converging in
$C(\overline\Om)\times C^1(\overline\Om)$.
Using the preceding estimates and the convergence at $t_j$,
we identify its limit as $(u_\infty,v_\infty)$.
Since the choice of $s_j$ was arbitrary, we conclude
\eqref{eq:compact-spatial}.
\end{proof}

Lemma~\ref{lem:subcritical-compactness} reduces convergence
to identifying the subsequential limits of $(u,v)$.
We use finite dissipation, the exponential decay of $w$, and radial
uniqueness to show that every such limit is
$(\overline{u_0},\overline{u_0})$.

\begin{lemma}\label{lem:subcritical-convergence}
We have $u(t)\to\overline{u_0}$ in $C(\overline\Om)$
and $v(t)\to\overline{u_0}$ in $C^1(\overline\Om)$
as $t\to\infty$.
More precisely,
\[
\|u(t)-\overline{u_0}\|_\infty
+\|v(t)-\overline{u_0}\|_{1,\infty}\to0
\qquad\text{as }t\to\infty.
\]
\end{lemma}

\begin{proof}
We fix a subsequence along which $(u,v)$ converges to
$(u_\infty,v_\infty)$ uniformly on $[t_j,t_j+1]$
in the sense of \eqref{eq:compact-spatial}.
For any $\varphi\in C^1(\overline\Om)$, we test
the $v$-equation of~\eqref{eq:model} against $\varphi$,
integrate by parts, and then integrate over
$[t_j,t_j+1]$ to obtain
\[
\begin{aligned}
\int_{t_j}^{t_j+1}\int_\Om
(\nabla v\cdot\nabla\varphi+v\varphi-u\varphi)\dd x\dd\tau&=-\int_{t_j}^{t_j+1}\int_\Om v_t\varphi\dd x\dd\tau\\
&=-\int_\Om\bigl(v(t_j+1)-v(t_j)\bigr)\varphi.
\end{aligned}
\]
Using the Cauchy--Schwarz inequality and the uniform
convergence on $[t_j,t_j+1]$ from
Lemma~\ref{lem:subcritical-compactness}, we deduce that
the right-hand side tends to zero: as $j\to \infty$,
\[
\left|\int_\Om\bigl(v(t_j+1)-v(t_j)\bigr)\varphi\right|
\le\|v(t_j+1)-v(t_j)\|_2\|\varphi\|_2\to0.
\]
Passing to the limit and using the uniform convergence in
Lemma~\ref{lem:subcritical-compactness}, we obtain
\[
\int_\Om
(\nabla v_\infty\cdot\nabla\varphi
+v_\infty\varphi-u_\infty\varphi)=0.
\]
Since $\varphi\in C^1(\overline\Om)$ is arbitrary,
we conclude that $v_\infty$ is a weak solution of
\begin{equation}\label{eq:stationary-limit-signal}
\begin{cases}
-\Delta v_\infty+v_\infty=u_\infty & \text{in }\Om,\\
\partial_\nu v_\infty=0 & \text{on }\partial\Om.
\end{cases}
\end{equation}
We next identify the limiting flux.
For any $\boldsymbol\varphi\in C_c^\infty(\Om;\mathbb R^2)$,
we combine H\"older's inequality with
\eqref{eq:vanishing-flux} to derive
\[
\begin{aligned}
&\left|\int_{t_j}^{t_j+1}\int_\Om
(\nabla u-\chi u\nabla v-\xi u\nabla w)
\cdot\boldsymbol\varphi\dd x\dd\tau\right|\\
\le&\|\boldsymbol\varphi\|_\infty
\int_{t_j}^{t_j+1}
\|\nabla u-\chi u\nabla v-\xi u\nabla w\|_1\dd\tau
\to0.
\end{aligned}
\]
We integrate the term involving $\nabla u$ by parts.
Passing to the limit and using the uniform convergence in
Lemma~\ref{lem:subcritical-compactness} and the decay estimate
for $\nabla w$ in \eqref{eq:preliminary-tissue-decay},
we arrive at
\[
-\int_\Om u_\infty\nabla\cdot\boldsymbol\varphi
-\chi\int_\Om u_\infty\nabla v_\infty
\cdot\boldsymbol\varphi=0.
\]
As $\boldsymbol\varphi$ is arbitrary, we conclude that
\begin{equation}\label{eq:stationary-limit-flux}
\nabla u_\infty=\chi u_\infty\nabla v_\infty
\quad\text{in distributions}.
\end{equation}
To identify the limit $(u_\infty,v_\infty)$, we rewrite
\eqref{eq:stationary-limit-flux} as
\[
\nabla(u_\infty e^{-\chi v_\infty})
=e^{-\chi v_\infty}
(\nabla u_\infty-\chi u_\infty\nabla v_\infty)=0.
\]
We therefore conclude that $u_\infty e^{-\chi v_\infty}$
is constant.
Using $\int_\Om u_\infty=\|u_0\|_1$, we deduce
\begin{equation}\label{eq:stationary-gibbs-relation}
u_\infty=\frac{\|u_0\|_1e^{\chi v_\infty}}
{\int_\Om e^{\chi v_\infty}}>0.
\end{equation}
Substituting \eqref{eq:stationary-gibbs-relation} into
\eqref{eq:stationary-limit-signal} and applying elliptic
regularity, we arrive at the classical equation
\begin{equation}\label{eq:stationary-mean-field}
\begin{cases}
\Delta(\chi v_\infty)-\chi v_\infty+
\frac{\chi\|u_0\|_1e^{\chi v_\infty}}
{\int_\Om e^{\chi v_\infty}}=0 & \text{in }\Om,\\
\partial_\nu v_\infty=0 & \text{on }\partial\Om.
\end{cases}
\end{equation}
In view of the radial symmetry of $v_\infty$ and
the condition $\chi\|u_0\|_1<8\pi$, we apply
the uniqueness theorem \cite[Theorem 1.1]{WangWangYang}
to \eqref{eq:stationary-mean-field}
and conclude that $v_\infty$ is constant.
Combining this fact with \eqref{eq:stationary-gibbs-relation}
and \eqref{eq:stationary-limit-signal}, we discover that
$u_\infty=v_\infty=\overline{u_0}$.
Since every sub-sequential limit has this value,
we conclude  the asserted uniform convergence from
Lemma~\ref{lem:subcritical-compactness}.
\end{proof}

Lemma~\ref{lem:subcritical-convergence} shows that $u$ and $v$
converge to $\overline{u_0}$.
We now prove that this convergence is exponential in
$L^2(\Om)\times H^1(\Om)$ by a quadratic energy estimate.

\begin{lemma}\label{lem:l2-stabilization}
There exist $C_{39}>0,\sigma_0>0$ such that
\[
\|u(t)-\overline{u_0}\|_2
+\|v(t)-\overline{u_0}\|_{1,2}
\le C_{39}e^{-\sigma_0t},\qquad\forall t\ge0.
\]
\end{lemma}

\begin{proof}
We test the $u$- and $v$-equations of~\eqref{eq:model} by
$u-\overline{u_0}$ and $-\chi\overline{u_0}\Delta v$,
respectively, and integrate by parts using the Neumann boundary
conditions to obtain
\begin{equation}\label{eq:quadratic-energy-identity}
\begin{aligned}
&\frac12\frac d{dt}
\left(\int_\Om(u-\overline{u_0})^2
+\chi\overline{u_0}\int_\Om|\nabla v|^2\right)\\
&+\int_\Om|\nabla u|^2
+\chi\overline{u_0}
\int_\Om|\Delta v|^2+\chi\overline{u_0}\int_\Om|\nabla v|^2\\
=&2\chi\overline{u_0}\int_\Om\nabla u\cdot\nabla v
+\chi\int_\Om(u-\overline{u_0})\nabla u\cdot\nabla v
+\xi\int_\Om u\nabla u\cdot\nabla w.
\end{aligned}
\end{equation}
We begin to estimate the first term on the right.
By the assumption $\chi\|u_0\|_1<8\pi$ and  the radial Neumann eigenvalue formula, we have
\[
\chi\overline{u_0}
=\frac{\chi\|u_0\|_1}{\pi R^2}
<\frac8{R^2}
<\lambda_1=\frac{j_{1,1}^2}{R^2}\approx \frac{14.6819}{R^2},
\]
where $j_{1,1}$ is the first positive zero of $J_1$;
see \cite{Watson}.
Applying the radial Poincar\'e inequality, we obtain
\begin{equation}\label{eq:cell-poincare}
\int_\Om|\nabla u|^2
\ge\lambda_1\int_\Om(u-\overline{u_0})^2.
\end{equation}
For $v$, we apply the same inequality to $v-\overline v$
to get
\begin{equation*}
\left(\int_\Om(v-\overline v)^2\right)^{1/2}
\le\lambda_1^{-1/2}
\left(\int_\Om|\nabla v|^2\right)^{1/2}.
\end{equation*}
We then  integrate by parts and apply the Cauchy--Schwarz
inequality  to deduce
\[
\|\nabla v\|_2^2=-\int_\Om(v-\overline v)\Delta v\le\|v-\overline v\|_2\|\Delta v\|_2\le\lambda_1^{-1/2}\|\nabla v\|_2\|\Delta v\|_2.
\]
Consequently, we have
\[
\int_\Om|\Delta v|^2
\ge\lambda_1\int_\Om|\nabla v|^2.
\]
In light of this estimate, we apply Young's inequality to obtain
\begin{equation}\label{eq:linear-coupling-estimate}
\begin{aligned}
2\chi\overline{u_0}
\left|\int_\Om\nabla u\cdot\nabla v\right|
&\le\frac{\chi\overline{u_0}}{\lambda_1}
\int_\Om|\nabla u|^2
+\chi\overline{u_0}\lambda_1\int_\Om|\nabla v|^2\\
&\le\frac{\chi\overline{u_0}}{\lambda_1}
\int_\Om|\nabla u|^2
+\chi\overline{u_0}\int_\Om|\Delta v|^2.
\end{aligned}
\end{equation}
We turn to the two remaining coupling terms in
\eqref{eq:quadratic-energy-identity}.
Fix $0<\epsilon<\frac14(1-\chi\overline{u_0}/\lambda_1)$.
Using the uniform convergence in Lemma~\ref{lem:subcritical-convergence},
we choose $T>0$ sufficiently large and consider $t\ge T$
in the following estimates.
Young's inequality then gives
\begin{equation}\label{eq:nonlinear-coupling-estimate}
\begin{aligned}
\chi\left|\int_\Om
(u-\overline{u_0})\nabla u\cdot\nabla v\right|
\le&\epsilon\int_\Om|\nabla u|^2
+\frac{\chi^2}{4\epsilon}
\|u-\overline{u_0}\|_\infty^2\int_\Om|\nabla v|^2\\
\le&\epsilon\int_\Om|\nabla u|^2
+\frac{\chi\overline{u_0}}4\int_\Om|\nabla v|^2.
\end{aligned}
\end{equation}
We also use the uniform bound for $u$ and the exponential decay estimate
\eqref{eq:preliminary-tissue-decay} to choose
a constant $C_{40}>0$ such that
\begin{equation}\label{eq:tissue-coupling-estimate}
\begin{aligned}
\xi\left|\int_\Om u\nabla u\cdot\nabla w\right|
\le&\epsilon\int_\Om|\nabla u|^2
+\frac{\xi^2}{4\epsilon}\|u\|_\infty^2
\int_\Om|\nabla w|^2\\
\le&\epsilon\int_\Om|\nabla u|^2
+C_{40}e^{-(\eta^0-\eta)t}.
\end{aligned}
\end{equation}
Substituting \eqref{eq:linear-coupling-estimate},
\eqref{eq:nonlinear-coupling-estimate}, and
\eqref{eq:tissue-coupling-estimate} into
\eqref{eq:quadratic-energy-identity} and rearranging,
we arrive at
\[
\begin{aligned}
&\frac12\frac d{dt}
\left(\int_\Om(u-\overline{u_0})^2
+\chi\overline{u_0}\int_\Om|\nabla v|^2\right)\\
&+\left(1-\frac{\chi\overline{u_0}}{\lambda_1}-2\epsilon\right)
\int_\Om|\nabla u|^2
+\frac{3\chi\overline{u_0}}4\int_\Om|\nabla v|^2
\le C_{40}e^{-(\eta^0-\eta)t}.
\end{aligned}
\]
With our choice of $\epsilon$, we have
\[
2\left(1-\frac{\chi\overline{u_0}}{\lambda_1}-2\epsilon\right)
\ge1-\frac{\chi\overline{u_0}}{\lambda_1}>0.
\]
Thus, with
$\sigma=\min\{\lambda_1-\chi\overline{u_0},\  3/2\}>0$, we apply the radial Poincar\'e inequality \eqref{eq:cell-poincare}
to deduce
\begin{equation}\label{eq:stability-energy}
\begin{aligned}
&\frac d{dt}
\left(\int_\Om(u-\overline{u_0})^2
+\chi\overline{u_0}\int_\Om|\nabla v|^2\right)\\
&+\sigma\left(\int_\Om(u-\overline{u_0})^2
+\chi\overline{u_0}\int_\Om|\nabla v|^2\right)
\le 2C_{40}e^{-(\eta^0-\eta)t}, \ \  \forall t\geq T.
\end{aligned}
\end{equation}
We choose $0<\sigma_0<\min\{\sigma/2,\ (\eta^0-\eta)/2,\ 1\}$.
Solving the  Gronwall's inequality \eqref{eq:stability-energy}
on $[T,t]$ via integrating factor method and using the continuity of the energy on $[0,T]$,
we find a constant $C_{41}>0$ such that, for $t\ge0$,
\begin{equation}\label{eq:quadratic-energy-decay}
\int_\Om(u-\overline{u_0})^2
+\chi\overline{u_0}\int_\Om|\nabla v|^2
\le C_{41}e^{-2\sigma_0t}.
\end{equation}
To obtain the $L^2$ decay of $v-\overline{u_0}$,
we also estimate its spatial mean.
Integrating the $v$-equation of~\eqref{eq:model} over $\Om$ and  using $\|u\|_1=\|u_0\|_1$,
we have
\[
\frac d{dt}\int_\Om(v-u_0)
=-\int_\Om(v-u_0),
\]
which simply gives
\[
\int_\Om(v-u_0)
=e^{-t}\int_\Om(v_0-u_0).
\]
Recalling that $\sigma_0<1$, we combine this identity with
\eqref{eq:quadratic-energy-decay} and the radial Poincar\'e
inequality to choose a constant $C_{42}>0$ such that, for $t\ge0$,
\begin{equation}\label{eq:signal-l2-decay}
\|v-\overline{u_0}\|_2\le\lambda_1^{-1/2}\|\nabla v\|_2
+\frac{e^{-t}}{|\Om|^{1/2}}
\left|\int_\Om(v_0-u_0)\right|\le C_{42}e^{-\sigma_0t}.
\end{equation}
Combining \eqref{eq:quadratic-energy-decay} with
\eqref{eq:signal-l2-decay}, we complete the proof.
\end{proof}

 Now,  one may further employ the standard parabolic Schauder theory to derive uniform bounds  for solutions in $C^{2+\alpha, 1+\alpha/2}(\overline{\Om}\times [t, t+1])$, and then, apply Gagliardo--Nirenberg interpolation inequalities  between Sobolev norms
 to lift the $L^2$-convergence to $L^\infty$-convergence, as commonly done in literature.  Herein, we wish to achieve the convergence lifting  based only on  H\"older bounds, to this end, we  present interpolation inequalities involving  Sobolev norms and H\"older norms. 
 \begin{lemma}\label{GN-Sob-hol}
 Let $k\in \mathbb{N}, \alpha\in (0,1)$ and $\Om\subset \mathbb{R}^n$ be a bounded and smooth domain. Then there exists a geometric constant $\Gamma=\Gamma(\Om, n, k, \alpha)>0$ such that
 \begin{equation}\label{GN-Soblev-holder}
\|D^kf\|_{L^\infty(\Omega)}\leq \Gamma\|f\|_{C^{k+\alpha}(\overline{\Om})}^\frac{2k+n}{2k+2\alpha+n}\|f\|_{L^2(\Omega)}^\frac{2\alpha}{2k+2\alpha+n}, \ \  \forall f\in C^{k+\alpha}(\overline{\Om}).
\end{equation}
\end{lemma}
\begin{proof}
To avoid triviality, we   proceed with  $\|D^k f\|_\infty>0$  and define
\begin{equation}\label{gM-def}
g=\frac{f}{\|f\|_{C^{k+\alpha}(\overline{\Om})}}, \quad \quad M=\|D^kg\|_{L^\infty(\Om)},
\end{equation}
so that
\[
\|g\|_{C^{k+\alpha}(\overline{\Om})}=1, \quad \quad  \|g\|_{L^2(\Om)}=\frac{\|f\|_{L^2(\Om)}}{\|f\|_{C^{k+\alpha}(\overline{\Om})}}.
\]
Then there exist $x_0\in\overline{\Om}$ and  $\beta \in \mathbb{N}^n$ with $|\beta|=k$ such that $|D^\beta g(x_0)|=M$.   Since the bounded and smooth domain $\Om$ satisfies  a uniform interior cone condition: there are geometric constants  $R_0, a>0$ and a fixed cone $K_0$ with vertex $0$ and nonempty interior such that for every $\rho\in (0, R_0]$,
\[
K:=x_0+\rho K_0\subset B(x_0, \rho) \cap \Om, \ \ \  |K|\ge a\rho^n.
\]
Let $P$ be the Taylor polynomial of $g$  at $x_0$ of degree $k$.  Since $\|D^\beta g\|_{C^\alpha(\overline{\Om})}\le 1$, there exists $\gamma=\gamma(R_0, K_0, k, \alpha)>0$ such that
\[
|g(y)-P(y)|\le \gamma \rho^{k+\alpha}, \  \ \ \forall y\in K.
\]
Thus,
\begin{equation}\label{Taylor-rem}
\|g-P\|_{L^2(K)}\le \gamma \rho^{k+\alpha}|K|^{1/2}\le \gamma |B(0,1)|^{1/2} \rho^{k+\alpha+n/2}.
\end{equation}
On the other hand, by finite-dimensional norm equivalence on the fixed cone $K_0$, there exists $b=b(K_0,k, n)>0$ such that
\[
\|P\|_{L^2(K_0)}\ge b|D^kP(0)|,
\]
and so, by changing variables with substitution  $Q(y)=P(x_0+\rho y)$, we find
\begin{align*}
\|P\|_{L^2(K)}=\rho^{n/2}\|Q\|_{L^2(K_0)}&\ge b\rho^{n/2}|D^kQ(0)|\\
&=b\rho^{k+n/2}|D^kP(x_0)|= b\rho^{k+ n/2}M.
\end{align*}
Next, in view of \eqref{Taylor-rem}, we choose small  $\epsilon>0$ with  $\rho:=(\epsilon M)^{1/\alpha}<R_0$ to satisfy
\[
\|g-P\|_{L^2(K)}\le \frac{1}{2} b\rho^{k+ n/2}M.
\]
Therefore,
\[
\|g\|_{L^2(K)}\ge \frac12 b\rho^{k+ n/2}M=\frac b2 \epsilon^\frac{k+n/2}{\alpha} M^{1+\frac{k+n/2}{\alpha}}.
\]
From this and the fact $K\subset \Omega$, we deduce from  \eqref{gM-def} that
\begin{align*}
\frac{\|f\|_{L^2(\Om)}}{\|f\|_{C^{k+\alpha}(\overline{\Om})}}=\|g\|_{L^2(\Om)}&\ge \frac b2 \epsilon^\frac{k+n/2}{\alpha} M^{1+\frac{k+n/2}{\alpha}}:=cM^{1+\frac{k+n/2}{\alpha}}\\
&=c\left(\frac{\|D^kf\|_{L^\infty(\Om)}}{\|f\|_{C^{k+\alpha}(\overline{\Om})}}\right)^{1+\frac{k+n/2}{\alpha}},
\end{align*}
which upon setting $\Gamma=c^{-\frac{\alpha}{k+\alpha+n/2}}$ implies \eqref{GN-Soblev-holder}.
\end{proof}
Next, we employ the uniform H\"older bounds  and the interpolation inequality  in Lemma ~\ref{GN-Sob-hol}  to strengthen the preceding
decay estimate to $L^\infty(\Om)\times W^{1,\infty}(\Om)$.
\begin{lemma}\label{lem:exponential-stabilization}
There exist $C_{43},\sigma>0$ such that
\[
\|u(t)-\overline{u_0}\|_\infty
+\|v(t)-\overline{u_0}\|_{1,\infty}
\le C_{43}e^{-\sigma t},\qquad\forall t\ge0.
\]
\end{lemma}

\begin{proof}
Let $\alpha$ be the H\"older exponent in
\eqref{eq:uniform-holder-bounds}.
Then $u-\overline{u_0}$ and $v-\overline{u_0}$ are uniformly
bounded in $C^\alpha(\overline\Om)$ and
$C^{1+\alpha}(\overline\Om)$, respectively, for $t\ge1$.
From Lemma~\ref{lem:l2-stabilization}, we have
\[
\|u-\overline{u_0}\|_2+\|v-\overline{u_0}\|_2
\le C_{39}e^{-\sigma_0t},\qquad t\ge0.
\]
Applying spatial interpolation inequality \eqref{GN-Soblev-holder} with $n=2$ and using the uniform
H\"older bounds, we find constants $C_{44},C_{45}, C_{46}, C_{47}>0$
such that, for $t\ge1$,
\[
\|u-\overline{u_0}\|_\infty
\le C_{44}
\|u-\overline{u_0}\|_{C^\alpha(\overline\Om)}^{1/(1+\alpha)}
\|u-\overline{u_0}\|_2^{\alpha/(1+\alpha)}\le C_{45}
e^{-\alpha\sigma_0t/(1+\alpha)},
\]
and
\[
\|v-\overline{u_0}\|_{1,\infty}
\le C_{46}
\|v-\overline{u_0}\|_{C^{1+\alpha}(\overline\Om)}^{2/(2+\alpha)}
\|v-\overline{u_0}\|_2^{\alpha/(2+\alpha)}\le C_{47}
e^{-\alpha\sigma_0t/(2+\alpha)}.
\]
Setting  $\sigma=\alpha\sigma_0/(2+\alpha)>0$,  we add the two estimates to deduce
\[
\|u-\overline{u_0}\|_\infty
+\|v-\overline{u_0}\|_{1,\infty}
\le \left(
C_{45}+C_{47}
\right)e^{-\sigma t},
\qquad t\ge1.
\]
By  the boundedness
of the solution on $[0,1]$, which gives
\[
\sup_{0\le t\le1}e^{\sigma t}
\left(
\|u-\overline{u_0}\|_\infty
+\|v-\overline{u_0}\|_{1,\infty}
\right)<\infty,
\]
we combine  the bounds on $[0,1]$ and $[1,\infty)$ to complete the proof.
\end{proof}

\begin{proof}[Proof of Theorem~\ref{thm:decay}]
By Lemma~\ref{lem:exponential-stabilization}, we have the asserted
exponential convergence of $u$ and $v$.
Integrating the decay estimate for $v$, we obtain
\[
\int_0^\infty\|v(\tau)-\overline{u_0}\|_{1,\infty}\dd\tau
\le\frac{C_{43}}{\sigma}.
\]
This readily entails
\[
\int_0^t v(x,\tau)\dd\tau
\ge \overline{u_0}t
-\int_0^\infty\|v(\tau)-\overline{u_0}\|_\infty\dd\tau
\ge \overline{u_0}t-\frac{C_{43}}{\sigma},
\]
and therefore,
\[
e^{\eta t-\int_0^t v(x,\tau)\dd\tau}
\le e^{C_{43}/\sigma}e^{-(\overline{u_0}-\eta)t}, \quad \forall t\ge 0.
\]
We now use the explicit formula \eqref{eq:w-exp} and the
gradient estimate \eqref{gradw-bdd} established in the proof of
Lemma~\ref{lem:bounded-bootstrap}.
Combining the preceding exponential bound with
$\int_0^t\|\nabla v(\tau)\|_\infty\dd\tau\le C_{43}/\sigma$,
we choose a constant $C_{48}>0$ such that, for all $t>0$,
\begin{align*}
\|w\|_{1,\infty}
\le & e^{C_{43}/\sigma}
\left(\|w_0\|_{1,\infty}
+\|w_0\|_\infty\int_0^t\|\nabla v(\tau)\|_\infty\dd\tau\right)
e^{-(\overline{u_0}-\eta)t}\\
\le & C_{48}e^{-(\overline{u_0}-\eta)t}.
\end{align*}
We have thus established \eqref{eq:subcritical-tissue-decay}
and completed the proof.
\end{proof}

\section{Finite-time blow-up}\label{sec:blowup}
In this section, we prove Theorem~\ref{thm:mass}.
We fix $\chi>0$, $\xi,\eta\ge0$, and $B>0$, and take
$w_0\in\mathcal W_B^0$.
Throughout this section, $(u,v,w)$ denotes the maximal radial
classical solution of~\eqref{eq:model} on $[0,T_m)$.
We derive an energy criterion for finite-time blow-up
and construct initial data satisfying it at every prescribed
mass above $8\pi/\chi$.

\begin{lemma}\label{lem:tissue}
For every $T_0>0$, there exists  a constant
$C_{49}(T_0,\|u_0\|_1,\|v_0\|_1)>0$ such that 
\begin{equation}\label{eq:tissue-gradient}
\|\nabla w(t)\|_\infty\le C_{49}B,
\qquad \forall t\in[0,\min\{T_0,T_m\}),\quad \forall w_0\in\mathcal W_B^0.
\end{equation}
\end{lemma}

\begin{proof}
By radiality and the divergence theorem, we find
\begin{equation*}
\int_{B_r}\Delta v(s)
=\int_{\partial B_r}\partial_\nu v(s)\dd S
=2\pi r v_r(r,s).
\end{equation*}
Then integrating the $v$-equation of~\eqref{eq:model} over
$B_r\times(0,t)$,  we obtain
\[
\begin{aligned}
2\pi r\partial_r\int_0^t v(r,s)\dd s
&=\int_0^t\int_{B_r}(v_s+v-u)\dd s\\
&=\int_{B_r}\left(v(t)+\int_0^t v(s)\dd s\right)
-\int_{B_r}\left(v_0+\int_0^t u(s)\dd s\right).
\end{aligned}
\]
We bound each of the two integrals on the right by
$(1+t)\max\{\|u_0\|_1,\|v_0\|_1\}$.
Their nonnegativity allows us to bound the absolute value
of their difference by the same quantity.
After dividing by $2\pi r$, we arrive at
\begin{equation}\label{eq:q-gradient}
\left|\partial_r\int_0^t v(r,s)\dd s\right|
\le\frac{(1+t)\max\{\|u_0\|_1,\|v_0\|_1\}}{2\pi r}.
\end{equation}
We now turn to the estimate for $w_r$.
Differentiating the $w$-equation of~\eqref{eq:model}
with respect to $r$, we have
\[
(w_r)_t=(\eta-v-2\eta w)w_r-wv_r.
\]
Solving this linear equation and using  the $w$-formula in \eqref{w-gradw-sol},
we conclude 
\begin{equation}\label{eq:w-radial-variation}
w_r(r,t)=e^{\eta t-\int_0^t(v(r,s)+2\eta w(r,s))\dd s}w_0'(r)-w(r,t)\int_0^t
e^{-\eta\int_s^t w(r,\tau)\dd\tau}v_r(r,s)\dd s.
\end{equation}
To estimate the integral in \eqref{eq:w-radial-variation},
we differentiate the exponential weight to find
\[
\partial_s e^{-\eta\int_s^t w(r,\tau)\dd\tau}
=\eta w(r,s)e^{-\eta\int_s^t w(r,\tau)\dd\tau}.
\]
This immediately gives rise to
\begin{equation}\label{eq:radial-weight-integral}
\eta\int_0^t
w(r,s)e^{-\eta\int_s^t w(r,\tau)\dd\tau}\dd s=\left[
e^{-\eta\int_s^t w(r,\tau)\dd\tau}
\right]_{s=0}^{s=t}=1-e^{-\eta\int_0^t w(r,s)\dd s}\le1.
\end{equation}
Writing
\[
v_r(r,s)=\partial_s\left(
\partial_r\int_0^s v(r,\tau)\dd\tau\right),
\]
we integrate by parts with respect to $s$ and deduce
\begin{equation}\label{eq:weighted-vr-identity}
\begin{aligned}
&\int_0^t e^{-\eta\int_s^t w(r,\tau)\dd\tau}v_r(r,s)\dd s\\
=&\left[
e^{-\eta\int_s^t w(r,\tau)\dd\tau}
\partial_r\int_0^s v(r,\tau)\dd\tau
\right]_{s=0}^{s=t}\\
&-\eta\int_0^t
w(r,s)e^{-\eta\int_s^t w(r,\tau)\dd\tau}
\left(\partial_r\int_0^s v(r,\tau)\dd\tau\right)\dd s\\
=&\partial_r\int_0^t v(r,s)\dd s-\eta\int_0^t
w(r,s)e^{-\eta\int_s^t w(r,\tau)\dd\tau}
\left(\partial_r\int_0^s v(r,\tau)\dd\tau\right)\dd s.
\end{aligned}
\end{equation}
We substitute \eqref{eq:weighted-vr-identity} into
\eqref{eq:w-radial-variation} and apply
\eqref{eq:radial-weight-integral} to control the time integral.
Combining this with $w(r,t)\le w_0(r)e^{\eta t}$
and \eqref{eq:q-gradient}, we infer 
\begin{equation}\label{eq:tissue-radial-gradient-bound}
\begin{aligned}
|w_r(r,t)|
&\le e^{\eta t}|w_0'(r)|
+2w(r,t)\sup_{0\le s\le t}
\left|\partial_r\int_0^s v(r,\tau)\dd\tau\right|\\
&\le e^{\eta t}|w_0'(r)|
+\frac{w_0(r)e^{\eta t}(1+t)}{\pi r}
\max\{\|u_0\|_1,\|v_0\|_1\}.
\end{aligned}
\end{equation}
Finally, since $w_0\in\mathcal W_B^0$ and $w_0(0)=0$, we have
\begin{equation}\label{eq:initial-tissue-radial-bounds}
|w_0'(r)|\le B,
\qquad
0\le w_0(r)\le Br,
\qquad\forall r\in[0,R].
\end{equation}
Combining \eqref{eq:tissue-radial-gradient-bound} with
\eqref{eq:initial-tissue-radial-bounds} and using $t\le T_0$, we obtain
\[
\|\nabla w(t)\|_\infty
\le Be^{\eta T_0}\left(
1+\frac{1+T_0}{\pi}\max\{\|u_0\|_1,\|v_0\|_1\}
\right)
\le C_{49}B.\qedhere
\]
\end{proof}

To control the negative part of the energy $\mathcal{F}$, we first establish
a lower bound for the logarithmically weighted entropy of $u$.

\begin{lemma}\label{lem:radial-log-entropy}
For every $q\in(1,2)$ and $\mu>0$, there exists a constant
$C_{50}(q,\mu)\ge e$ such that, whenever
$\chi\max\{\|u\|_1,\|v\|_{1,q}\}\le\mu$, we have
\begin{equation}\label{eq:radial-log-entropy}
2\pi\int_0^R ru(r)\ln\bigl(2\pi r^2u(r)\bigr)\dd r
+\frac{\mu}{\chi}\ln(1+\ln A)
\ge-\left(1+\frac{8\pi}{\chi}+\frac4e\right),
\end{equation}
where $A$ is given by 
\[
A:=C_{50}\left[
1+\|v_t\|_2
+\left(\int_\Om u|\nabla(\ln u-\chi v)|^2\right)^{1/2}
\right].
\]
\end{lemma}

\begin{proof}
Integrating the $v$-equation of~\eqref{eq:model} over $B_r$,
we obtain
\begin{equation}\label{eq:radial-signal-identity}
2\pi rv_r(r)
=\int_{\partial B_r}\partial_\nu v\dd S
=\int_{B_r}\Delta v
=\int_{B_r}(v_t+v)-\int_{B_r}u.
\end{equation}
For $u$, we compute in polar coordinates
\[
\frac d{dr}\int_{B_r}u
=2\pi\frac d{dr}\int_0^r\rho u(\rho)\dd\rho
=2\pi ru(r).
\]
By the fundamental theorem of calculus, we then deduce
\begin{equation}\label{eq:radial-mass-square}
\begin{aligned}
\int_0^r\rho u(\rho)\int_{B_\rho}u\,\dd x\,\dd\rho
&=\frac1{2\pi}\int_0^r
\int_{B_\rho}u\,\dd x\,
\frac d{d\rho}\left(\int_{B_\rho}u\,\dd x\right)\dd\rho\\
&=\frac1{4\pi}\left(\int_{B_r}u\,\dd x\right)^2.
\end{aligned}
\end{equation}
Integrating by parts and applying
\eqref{eq:radial-signal-identity} and
\eqref{eq:radial-mass-square}, we find
\begin{equation}\label{eq:radial-mass-flux}
\begin{aligned}
&\int_{B_r}u(x)\,x\cdot\nabla(\ln u-\chi v)(x)\dd x\\
=&\int_{B_r}x\cdot\nabla u
-\chi\int_{B_r}u(x)\,x\cdot\nabla v(x)\dd x\\
=&\int_{\partial B_r}u(x)\,x\cdot\nu\dd S
-2\int_{B_r}u
-2\pi\chi\int_0^r\rho^2u(\rho)v_r(\rho)\dd\rho\\
=&2\pi r^2u(r)-2\int_{B_r}u
-\chi\int_0^r\rho u(\rho)
\left(\int_{B_\rho}(v_t+v)-\int_{B_\rho}u\right)\dd\rho\\
=&2\pi r^2u(r)
-\left(2-\frac{\chi}{4\pi}\int_{B_r}u\right)\int_{B_r}u
-\frac{\chi}{2\pi}\int_{B_r}u(x)
\left(\int_{B_{|x|}}(v_t+v)\right)\dd x.
\end{aligned}
\end{equation}
To estimate the integral terms in \eqref{eq:radial-mass-flux},
we apply the Cauchy--Schwarz inequality to derive 
\begin{equation}\label{eq:radial-flux-cauchy}
\left|\int_{B_r}u(x)\,x\cdot\nabla(\ln u-\chi v)(x)\dd x\right|
\le r\|u\|_1^{1/2}
\left(\int_\Om u|\nabla(\ln u-\chi v)|^2\right)^{1/2}
\end{equation}
and
\begin{equation}\label{eq:radial-signal-cauchy}
\left|\int_{B_\rho}(v_t+v)\right|
\le\sqrt\pi\rho\bigl(\|v_t\|_2+\|v\|_2\bigr),
\qquad\forall\rho\in(0,r].
\end{equation}
We substitute \eqref{eq:radial-flux-cauchy} and
\eqref{eq:radial-signal-cauchy} into
\eqref{eq:radial-mass-flux} and use the embedding
$W^{1,q}(\Omega)\hookrightarrow L^2(\Omega)$
to control $\|v\|_2$.
We then choose $C_{50}>e$ sufficiently large to obtain
\begin{equation}\label{eq:radial-mass-error}
\begin{aligned}
&\left|2\pi r^2u(r)
-\left(2-\frac{\chi}{4\pi}\int_{B_r}u\right)
\int_{B_r}u\right|\\
\le&\frac{C_{50}r}{R}\left[
1+\|v_t\|_2
+\left(\int_\Om u|\nabla(\ln u-\chi v)|^2\right)^{1/2}
\right]
=\frac{Ar}{R},
\qquad\forall r\in(0,R].
\end{aligned}
\end{equation}
We now split the target integral
$2\pi\int_0^R ru(r)\ln(2\pi r^2u(r))\dd r$
at $r=R/A$.
On the outer interval $(R/A,R)$, the mass bound gives
$0<2\pi\int_{R/A}^R ru(r)\dd r\le\|u\|_1\le \mu/\chi$.
Applying Jensen's inequality to the convex function $f(z)=z\ln z$ for $z>0$ with respect to
$\dd r/(r\ln A)$ and using $s\ln s\ge-1/e$ and $A\ge e$,
we infer
\begin{equation}\label{eq:outer-log-entropy}
\begin{split}
&2\pi\int_{R/A}^R ru(r)\ln\bigl(2\pi r^2u(r)\bigr)\dd r\\
=&\int_{R/A}^R
\bigl(2\pi r^2u(r)\bigr)\ln\bigl(2\pi r^2u(r)\bigr)
\frac{\dd r}{r}\\
\ge&\left(2\pi\int_{R/A}^R ru(r)\dd r\right)
\ln\frac{2\pi\int_{R/A}^R ru(r)\dd r}{\ln A}\\
=&\left(2\pi\int_{R/A}^R ru(r)\dd r\right)
\ln\left(2\pi\int_{R/A}^R ru(r)\dd r\right)
-\left(2\pi\int_{R/A}^R ru(r)\dd r\right)\ln\ln A\\
\ge&-\frac1e-\frac{\mu}{\chi}\ln(1+\ln A).
\end{split}
\end{equation}
On the inner interval $(0,R/A)$, we define $E\subset(0,R/A)$ as the set of radii satisfying
\begin{equation}\label{E-def}
0<\int_{B_r}u<\frac{8\pi}{\chi},
\qquad
2\pi r^2u(r)\ge
\left(1-\frac{\chi}{8\pi}\int_{B_r}u\right)\int_{B_r}u.
\end{equation}
We estimate the integral over $E$ from below.
Under the substitution $\sigma=\int_{B_r}u$, with
$\dd\sigma=2\pi ru(r)\dd r$, we denote the image of $E$ under $\sigma$
by $S\subset(0,8\pi/\chi)$.
Using the monotonicity of the  logarithm function, together with
$\ln\sigma\ge0$ for $\sigma\ge1$ and
$\ln(1-\chi\sigma/(8\pi))\le0$
for $0<\sigma< 8\pi/\chi$, we deduce
\begin{equation}\label{eq:inner-log-entropy-E}
\begin{aligned}
&2\pi\int_E ru(r)\ln\bigl(2\pi r^2u(r)\bigr)\dd r\\
\ge&2\pi\int_E ru(r)
\ln\left(
\left(1-\frac{\chi}{8\pi}\int_{B_r}u\right)\int_{B_r}u
\right)\dd r\\
=&\int_S
\ln\left(\sigma\left(1-\frac{\chi\sigma}{8\pi}\right)\right)
\dd\sigma\\
\ge&\int_0^{\min\{8\pi/\chi,1\}}\ln\sigma\dd\sigma
+\int_0^{8\pi/\chi}
\ln\left(1-\frac{\chi\sigma}{8\pi}\right)\dd\sigma\\
\ge&\int_0^1\ln\sigma\dd\sigma
+\frac{8\pi}{\chi}\int_0^1\ln z\dd z\\
=&-1-\frac{8\pi}{\chi}.
\end{aligned}
\end{equation}
For $r\in(0,R/A)\setminus E$, we use
the definition of $E$ in \eqref{E-def} and then apply
\eqref{eq:radial-mass-error} to derive
\[
0<2\pi r^2u(r)
\le\left|2\pi r^2u(r)
-\left(2-\frac{\chi}{4\pi}\int_{B_r}u\right)\int_{B_r}u\right|
\le\frac{Ar}{R}.
\]
For these radii with $0<r\le R/(eA)$, we have seen that
$2\pi r^2u(r)\le Ar/R\le e^{-1}$.
Using the increasing monotonicity of $-a\ln a$ on $(0,e^{-1}]$,
we get 
\[
-\bigl(2\pi r^2u(r)\bigr)\ln\bigl(2\pi r^2u(r)\bigr)
\le-\frac{Ar}{R}\ln\frac{Ar}{R}
=\frac{Ar}{R}\ln\frac{R}{Ar}.
\]
For $R/(eA)<r<R/A$, we instead use
$a\ln a\ge-e^{-1}$ for $a>0$.
Integrating these estimates and enlarging the integration
domains, we obtain
\begin{equation}\label{eq:inner-log-entropy-complement}
\begin{aligned}
&2\pi\int_{(0,R/A)\setminus E}
ru(r)\ln\bigl(2\pi r^2u(r)\bigr)\dd r\\
\ge&-\int_0^{R/(eA)}
\frac{Ar}{R}\ln\frac{R}{Ar}\,\frac{\dd r}{r}
-\int_{R/(eA)}^{R/A}\frac1e\,\frac{\dd r}{r}\\
=&-\frac3e.
\end{aligned}
\end{equation}
Combining \eqref{eq:outer-log-entropy},
\eqref{eq:inner-log-entropy-E}, and
\eqref{eq:inner-log-entropy-complement}, we conclude
\eqref{eq:radial-log-entropy}.
\end{proof}

We now combine this entropy bound with a weighted estimate
for $\nabla v$ to control the negative part of
$\mathcal F+\xi\int_\Om uw$, where $\mathcal F$ is defined
by  \eqref{eq:full-energy}.

\begin{lemma}\label{thm:spatial}
Let $q\in(1,2)$ and $\mu>0$, and assume that $\chi \max\left\{\|u\|_1,\|v\|_{1,q}\right\}\le\mu$.
Then there exists $C_{51}(q,\mu)>0$ such that
\begin{equation}\label{eq:spatial-energy}
-\mathcal F(u,v,w)-\xi\int_\Om uw
\le C_{51}\left(
1+\chi\int_\Omega v_t^2
+\int_\Om u|\nabla(\ln u-\chi v)|^2
\right)^{\frac34}.
\end{equation}
\end{lemma}

\begin{proof}
Throughout the proof, we write $r=|x|$.
We first estimate a weighted $L^2$ norm of $\nabla v$.
By polar integration, we have, for $q\in (1,2)$,
\[
\int_{R/2}^R
\bigl(|v(\rho)|^q+|v_r(\rho)|^q\bigr)\dd\rho
\le\frac2R\int_{R/2}^R
\rho\bigl(|v(\rho)|^q+|v_r(\rho)|^q\bigr)\dd\rho
\le\frac{\|v\|_{1,q}^q}{\pi R}
\le\frac{\mu^q}{\pi R\chi^q}.
\]
We therefore control $v$ in $W^{1,q}(R/2,R)$.
Applying the embedding
$W^{1,q}(R/2,R)\hookrightarrow C([R/2,R])$,
we obtain a uniform bound for $v(R)$.
By H\"older's inequality, we deduce
\[
\begin{aligned}
|v(r)-v(R)|
&\le\int_r^R
\rho^{1/q}|v_r(\rho)|\,\rho^{-1/q}\dd\rho\\
&\le
\left(\int_r^R\rho |v_r(\rho)|^q\right)^\frac{1}{q}\left(\int_r^R \rho^{-\frac{1}{q-1}}\right)^\frac{q-1}{q}\\
&\le \frac{\mu}{\chi(2\pi)^{1/q}}
\left(\frac{q-1}{2-q}\right)^{(q-1)/q}
r^{-(2-q)/q}.
\end{aligned}
\]
Combining these bounds, we find $C_{52}(q,\mu)\ge1$ such that
\[
|v(r)|\le |v(R)|+|v(r)-v(R)|
\le C_{52}\bigl(1+r^{-(2-q)/q}\bigr),
\qquad\forall r\in(0,R].
\]
As $q>1$, this  enables us to find  a constant $C_{53}(q,\mu)>0$
such that
\begin{equation}\label{eq:weighted-v-bounds}
\begin{aligned}
&\left\|r^2v\ln^2\frac rR\right\|_\infty
+\int_\Om
\left(1+\left|\ln\frac rR\right|\right)^2v^2\\
=&\sup_{0<r\le R}r^2|v(r)|\ln^2\frac rR
+2\pi\int_0^R
r\left(1+\left|\ln\frac rR\right|\right)^2v^2(r)\dd r\\
\le& 4\pi C^2_{52}\left[
\sup_{0<r\le R}
\left(r^2+r^{3-2/q}\right)\ln^2\frac rR
+\int_0^R
\left(r+r^{3-4/q}\right)
\left(1+\left|\ln\frac rR\right|\right)^2d r
\right]\\
\le& C_{53}.
\end{aligned}
\end{equation}
We next test $-\Delta v+v=u-v_t$ by
$r^2v\ln^2(r/R)$ and integrate by parts to obtain
\begin{equation}\label{eq:weighted-elliptic-testing}
\begin{aligned}
&\int_\Om r^2|\nabla v|^2\ln^2\frac rR
+\int_\Om r^2v^2\ln^2\frac rR\\
=&\int_\Om r^2uv\ln^2\frac rR
-\int_\Om r^2v_tv\ln^2\frac rR-\int_\Om
v\nabla\left(r^2\ln^2\frac rR\right)\cdot\nabla v.
\end{aligned}
\end{equation}
To estimate the last term, we note
\[
\left|\nabla\left(r^2\ln^2(r/R)\right)\right|^2=4\left(1+\ln\frac rR\right)^2r^2\ln^2(r/R),
\]
and then, we apply Young's inequality to see that
\begin{equation}\label{eq:weighted-cross-term}
\begin{aligned}
\left|\int_\Om
v\nabla\left(r^2\ln^2\frac rR\right)\cdot\nabla v\right|
\le&\frac12\int_\Om
r^2|\nabla v|^2\ln^2\frac rR
+\frac12\int_\Om
\frac{\left|\nabla\left(r^2\ln^2(r/R)\right)\right|^2}
{r^2\ln^2(r/R)}v^2\\
=&\frac12\int_\Om
r^2|\nabla v|^2\ln^2\frac rR
+2\int_\Om\left(1+\ln\frac rR\right)^2v^2.
\end{aligned}
\end{equation}
Substituting \eqref{eq:weighted-cross-term} into
\eqref{eq:weighted-elliptic-testing} and using
H\"older's inequality and \eqref{eq:weighted-v-bounds},
we find a constant $C_{54}(q,\mu)>0$ such that
\[
\begin{aligned}
\int_\Om r^2|\nabla v|^2\ln^2\frac rR
&\le
\|u\|_1\left\|r^2v\ln^2\frac rR\right\|_\infty
+\|v_t\|_2
\left\|r^2v\ln^2\frac rR\right\|_2\\
&\quad+\frac12\int_\Om
r^2|\nabla v|^2\ln^2\frac rR
+2\int_\Om\left(1+\ln\frac rR\right)^2v^2\\
&\le
\frac12\int_\Om r^2|\nabla v|^2\ln^2\frac rR
+C_{54}\bigl(1+\|v_t\|_2\bigr).
\end{aligned}
\]
Therefore, absorbing the gradient term and using radiality, we end up with
\begin{equation}\label{eq:weighted-signal-energy}
\int_\Om \left|x\cdot\nabla v\,\ln\frac rR\right|^2
\le 2C_{54}\bigl(1+\|v_t\|_2\bigr).
\end{equation}
By radial integration and integration by parts, we have
\begin{equation}\label{eq:weighted-radial-identity}
\begin{aligned}
\int_\Om\Delta v\,(x\cdot\nabla v)\ln\frac rR
&=2\pi\int_0^R
\left(v_{rr}+\frac1r v_r\right)
r^2v_r\ln\frac rR\,\dd r\\
&=\pi\int_0^R
\frac d{dr}\left(r^2v_r^2\right)\ln\frac rR\,\dd r\\
&=\pi\left[r^2v_r^2\ln\frac rR\right]_0^R
-\pi\int_0^R rv_r^2\dd r=-\frac12\int_\Om|\nabla v|^2.
\end{aligned}
\end{equation}
Testing $-\Delta v+v=u-v_t$ against $v$ and integrating
by parts, we also obtain
\begin{equation*}
\int_\Om uv
=\int_\Om(-\Delta v+v+v_t)v
=\int_\Om\bigl(|\nabla v|^2+v^2\bigr)+\int_\Om v_tv,
\end{equation*}
which along with \eqref{eq:weighted-radial-identity} turns the definition of
$\mathcal F$ in \eqref{eq:full-energy} into
\begin{equation}\label{eq: energy-identity-newform}
\mathcal F+\xi\int_\Om uw-\int_\Om u\ln u=\chi\int_\Om\Delta v\,(x\cdot\nabla v)\ln\frac rR-\frac\chi2\int_\Om v^2
-\chi\int_\Om v_tv.
\end{equation}
Integrating by parts gives
\[
\int_\Om x\cdot\nabla f \ln \frac rR=-\int_\Om f(1+2\ln \frac rR), \ \ \  \ f\in \{u, v^2\}.
\]
Then using
  $\Delta v=-u+v+v_t$, we further compute from \eqref{eq: energy-identity-newform} that
\begin{equation}\label{eq:weighted-energy-identity}
\begin{aligned}
&\mathcal F+\xi\int_\Om uw-\int_\Om u\ln u\\
={}&
-\chi\int_\Om u(x\cdot\nabla v)\ln\frac rR
+\frac\chi2\int_\Om x\cdot\nabla(v^2)\ln\frac rR
-\frac\chi2\int_\Om v^2\\
&+\chi\int_\Om v_t
\left(x\cdot\nabla v\,\ln\frac rR-v\right)\\
={}&
-\int_\Om x\cdot\nabla u\,\ln\frac rR+\int_\Om u\,x\cdot\nabla(\ln u-\chi v)\,\ln\frac rR\\
&-\chi\int_\Om\left(1+\ln\frac rR\right)v^2
+\chi\int_\Om v_t
\left(x\cdot\nabla v\,\ln\frac rR-v\right)\\
={}&2\int_\Om u\ln\frac rR
+\int_\Om u-\chi\int_\Om\left(1+\ln\frac rR\right)v^2\\
&+\int_\Om u\,x\cdot\nabla(\ln u-\chi v)\,\ln\frac rR
+\chi\int_\Om v_t
\left(x\cdot\nabla v\,\ln\left(\frac rR\right)-v\right).
\end{aligned}
\end{equation}
We further control the terms in $\mathcal{F}$. First, by polar integration, we have
\[
\int_\Om u\ln u+2\int_\Om u\ln\frac rR
=2\pi\int_0^R ru(r)\ln\bigl(2\pi r^2u(r)\bigr)\dd r
-\|u\|_1\ln(2\pi R^2).
\]
Second, applying the Cauchy--Schwarz inequality and the mass bound,
we obtain
\begin{equation*}
\begin{aligned}
\left|\int_\Om
u\,x\cdot\nabla(\ln u-\chi v)\,\ln\frac rR\right|
&\le\left\|r\ln\frac rR\right\|_\infty
\|u\|_1^{1/2}
\left(\int_\Om u|\nabla(\ln u-\chi v)|^2\right)^{1/2}\\
&\le \frac{R}{e}\|u\|_1^\frac12
\left(\int_\Om u|\nabla(\ln u-\chi v)|^2\right)^{1/2}.
\end{aligned}
\end{equation*}
Using \eqref{eq:weighted-signal-energy} and the $L^2$ bound
for $v$, we similarly obtain a constant $C_{55}(q,\mu)>0$ such that
\begin{equation*}
\begin{aligned}
\left|\int_\Om v_t
\left(x\cdot\nabla v\,\ln\left(\frac rR\right)-v\right)\right|
&\le\|v_t\|_2
\left(\left\|x\cdot\nabla v\,\ln\frac rR\right\|_2
+\|v\|_2\right)\\
&\le C_{55}\left(\|v_t\|_2+\|v_t\|_2^{3/2}\right).
\end{aligned}
\end{equation*}
Substituting the above  three estimates   into
\eqref{eq:weighted-energy-identity}, and then  applying the mass bound, \eqref{eq:weighted-v-bounds} and
Lemma~\ref{lem:radial-log-entropy}, we find a constant $C_{56}(q,\mu)>0$ such that
\[
-\mathcal F-\xi\int_\Om uw
\le C_{56}\biggl[
1+\ln(1+\ln A)+\|v_t\|_2^{3/2}
+\left(\int_\Om
u|\nabla(\ln u-\chi v)|^2\right)^{1/2}
\biggr].
\]
Using $\ln(1+\ln A)\le A$ and the definition of $A$
in Lemma~\ref{lem:radial-log-entropy}, we conclude
\eqref{eq:spatial-energy}.
\end{proof}

Combining the gradient bound of $w$ with the preceding energy
estimate, we derive a differential inequality for the energy
and obtain a criterion for finite-time blow-up.

\begin{lemma}\label{prop:criterion}
Let $q\in(1,2)$ and $\mu,T_0>0$, and assume that
$\chi\max\{\|u_0\|_1,\|v_0\|_{1,q}\}\le\mu$.
There exist $L(q,\mu)\ge1$ and
$C_{57}(q,\mu,T_0,B)>0$ such that, if
\begin{equation}\label{eq:thresholds}
\mathcal F(u_0,v_0,w_0)+\xi\int_\Om u_0w_0
<-\max\left\{(4LC_{57})^{3/4},
\left(\frac{12L}{T_0}\right)^3\right\},
\end{equation}
then the maximal existence time $T_m$ satisfies
\begin{equation}\label{eq:criterion-time}
T_m\le12L
\left(
-\mathcal F(u_0,v_0,w_0)-\xi\int_\Om u_0w_0
\right)^{-\frac{1}{3}}
<T_0,
\end{equation}
and the cell density blows up at $T_m$:
\begin{equation}\label{eq:criterion-blowup}
\limsup_{t\to T_m}\|u(t)\|_\infty=\infty.
\end{equation}
\end{lemma}

\begin{proof}
We test the $u$-equation of \eqref{eq:model} by
$\ln u-\chi v$ and integrate by parts:
\begin{equation}\label{eq:u-logarithmic-testing-identity}
\int_\Om u_t(\ln u-\chi v)
=-\int_\Om u|\nabla(\ln u-\chi v)|^2
+\xi\int_\Om u\nabla w\cdot\nabla(\ln u-\chi v).
\end{equation}
We add $\xi\int_\Om uw$ to the energy functional $\mathcal{F}$ defined in
\eqref{eq:full-energy} and differentiate with respect to $t$.
Substituting \eqref{eq:u-logarithmic-testing-identity},
integrating the gradient term by parts, and using
the $v$-equation  of \eqref{eq:model}, we conclude 
\begin{equation}\label{eq:modified-energy-identity}
\begin{aligned}
&\frac d{dt}\left(\mathcal F+\xi\int_\Om uw\right)\\
=&\int_\Om u_t(\ln u-\chi v)
+\chi\int_\Om\nabla v\cdot\nabla v_t
+\chi\int_\Om(v-u)v_t\\
=&-\int_\Om u|\nabla(\ln u-\chi v)|^2
+\xi\int_\Om u\nabla w\cdot\nabla(\ln u-\chi v)-\chi\int_\Om v_t^2.
\end{aligned}
\end{equation}
To estimate the term involving $\nabla w$ on the right-hand side of
\eqref{eq:modified-energy-identity}, we use the $L^\infty$-bound for
$\nabla w$ in \eqref{eq:tissue-gradient} and Young's inequality,
together with $\chi\max\{\|u_0\|_1,\|v_0\|_{1,q}\}\le\mu$, to obtain
\begin{equation*}
\begin{aligned}
\xi\int_\Om u\nabla w\cdot\nabla(\ln u-\chi v)
&\le\frac12\int_\Om u|\nabla(\ln u-\chi v)|^2
+\frac{\xi^2}{2}\int_\Om u|\nabla w|^2\\
&\le\frac12\int_\Om u|\nabla(\ln u-\chi v)|^2
+\frac{\xi^2\mu C_{49}^2B^2}{2\chi}.
\end{aligned}
\end{equation*}
Substituting this estimate  into
\eqref{eq:modified-energy-identity}, we find
\begin{equation}\label{eq:energy-error}
\frac d{dt}\left(\mathcal F+\xi\int_\Om uw\right)
\le-\frac12\left(
\chi\int_\Om v_t^2
+\int_\Om u|\nabla(\ln u-\chi v)|^2
\right)+\frac{\xi^2\mu C_{49}^2B^2}{2\chi}.
\end{equation}
We next use the variation-of-constants formula for $v$,
the Neumann semigroup estimates, and the assumption
$\chi \max\{\|u_0\|_1,\|v_0\|_{1,q}\}\le\mu$
to find, for every $q\in(1,2)$,
a constant $C_{58}(q)\ge1$ such that
$\max\{\|u(t)\|_1,\|v(t)\|_{1,q}\}\le C_{58}\mu/\chi$
for all $t\in[0,T_m)$.
Applying Lemma~\ref{thm:spatial} with $C_{58}\mu$ in place of
$\mu$, we obtain $L(q,\mu)\ge1$ such that
\[
\left(-\mathcal F-\xi\int_\Om uw\right)_+^{4/3}
\le L\left(
\chi\int_\Omega v_t^2
+\int_\Om u|\nabla(\ln u-\chi v)|^2
\right)+L,
\]
which, upon rearrangement, yields
\begin{equation}\label{eq:negative-dissipation-bound}
-\frac12\left(
\chi\int_\Om v_t^2
+\int_\Om u|\nabla(\ln u-\chi v)|^2
\right)
\le\frac12
-\frac{\left(-\mathcal F-\xi\int_\Om uw\right)_+^{4/3}}{2L}.
\end{equation}
Combining \eqref{eq:negative-dissipation-bound}
with \eqref{eq:energy-error}, we find a constant
$C_{57}(q,\mu,T_0,B)>0$ such that
\begin{equation}\label{eq:superlinear-energy}
\frac d{dt}\left(\mathcal F+\xi\int_\Om uw\right)
\le C_{57}
-\frac{1}{2L}\left(-\mathcal F-\xi\int_\Om uw\right)_+^{\frac43},
\ \ \forall t\in(0,\min\{T_0,T_m\}).
\end{equation}
Set
\[
y(t):=-\mathcal F(t)-\xi\int_\Om uw.
\]
By \eqref{eq:thresholds} and  \eqref{eq:superlinear-energy},
we infer
\begin{equation}\label{eq:y-differential-inequality}
\begin{cases}
y'(t)\ge\frac{y_+^{4/3}(t)}{2L}-C_{57},
\quad  \forall t\in(0,\min\{T_0,T_m\})\\[0.25cm]
y(0)>\max\left\{(4LC_{57})^{3/4},(12L/T_0)^3\right\}
\end{cases}
\end{equation}
Since $y(0)>(4LC_{57})^{3/4}$, a first time
$t_*\in(0,\min\{T_0,T_m\})$ at which $y$ reaches this
threshold from above would satisfy $y'(t_*)\le0$.
However, evaluating \eqref{eq:y-differential-inequality}
at $t=t_*$, we get 
\[
y'(t_*)\ge\frac{y^{4/3}(t_*)}{2L}-C_{57}
=C_{57}>0,
\]
a contradiction.
We therefore conclude that
\[
y(t)>(4LC_{57})^{3/4},\quad y'(t)\ge\frac{y^{4/3}(t)}{2L}-C_{57}\ge\frac{y^{4/3}(t)}{4L}>0,\ \forall t\in(0,\min\{T_0,T_m\}).
\]
Consequently,
\[
\frac d{dt}\left(y^{-1/3}(t)\right)
=-\frac13y^{-4/3}(t)y'(t)
\le-\frac{1}{3y^{4/3}(t)}
\cdot\frac{y^{4/3}(t)}{4L}
=-\frac1{12L}.
\]
Integrating from $0$ to $t$ yields
\begin{equation}\label{eq:y-integrated-inequality}
0<y^{-1/3}(t)
\le y^{-1/3}(0)-\frac{t}{12L},
\qquad\forall t\in[0,\min\{T_0,T_m\}).
\end{equation}
Notice also from \eqref{eq:y-differential-inequality} that $y(0)>(12L/T_0)^3$, and hence
$12Ly^{-1/3}(0)<T_0$.
Suppose to the contrary  that
$T_m>12Ly^{-1/3}(0)$.
Then $12Ly^{-1/3}(0)<\min\{T_0,T_m\}$,
so substituting $t=12Ly^{-1/3}(0)$ into
\eqref{eq:y-integrated-inequality} gives
\[
0<y^{-1/3}(t)
\le y^{-1/3}(0)
-\frac{12Ly^{-1/3}(0)}{12L}
=0,
\]
a contradiction.
Therefore,
$T_m\le12Ly^{-1/3}(0)<T_0$,
which proves \eqref{eq:criterion-time}.
Finally, suppose $\|u(t)\|_\infty$ remained bounded
as $t\to T_m$.
The Neumann semigroup estimates would also bound
$\|v(t)\|_{1,\infty}$.
These bounds contradict the extensibility criterion
\eqref{eq:local-extensibility} at the finite maximal
existence time.
Thus \eqref{eq:criterion-blowup} holds.
\end{proof}

To apply Lemma~\ref{prop:criterion}, we next construct
concentrated radial initial data with any prescribed mass
$m>8\pi/\chi$.
The following lemma shows that the corresponding modified
energy tends to $-\infty$ as $\varepsilon\to0+$.

\begin{lemma}\label{lem:initial-energy}
Let $m>8\pi/\chi$, and choose a nonnegative radial
function $\varphi\in C_c^\infty(B_1)$ satisfying
$\int_{B_1}\varphi=1$.
For $\varepsilon\in(0,\min\{R/2,\ 1,\ \sqrt{m\chi}\})$, define
\[
u_{0,\varepsilon}(x)
=\left(m-\frac{\varepsilon^2}\chi\right)
\varepsilon^{-2}\varphi(\frac x \varepsilon)
+\frac{\varepsilon^2}{\chi|\Om|},
\]
and let $v_{0,\varepsilon}$ solve
\begin{equation}\label{eq:concentrated-signal}
\begin{cases}
-\Delta v_{0,\varepsilon}+v_{0,\varepsilon}=u_{0,\varepsilon}
&\text{in }\Om,\\[0.2cm]
\partial_\nu v_{0,\varepsilon}=0
&\text{on }\partial\Om.
\end{cases}
\end{equation}
Then $\|u_{0,\varepsilon}\|_1=m$.
Moreover, as $\varepsilon\to0$,
\begin{equation}\label{eq:energy-expansion}
\mathcal F(u_{0,\varepsilon},v_{0,\varepsilon},w_0)
+\xi\int_\Om u_{0,\varepsilon}w_0
=-\frac{m( m\chi-8\pi)}{4\pi}
\ln(\frac 1\varepsilon)+O(1)\to-\infty.
\end{equation}
\end{lemma}

\begin{proof}
We begin with a semigroup representation of $v_{0,\varepsilon}$.
By the $L^\infty$ contractivity of the Neumann heat semigroup implied by the maximum principle,
we have
\[
\left\|e^{s(\Delta-1)}u_{0,\varepsilon}\right\|_\infty
=e^{-s}\left\|e^{s\Delta}u_{0,\varepsilon}\right\|_\infty
\le e^{-s}\|u_{0,\varepsilon}\|_\infty\to0,
\qquad\text{as }s\to\infty.
\]
Using the identity
$$\frac d{ds}\left(e^{s(\Delta-1)}u_{0,\varepsilon}\right)
=(\Delta-1)e^{s(\Delta-1)}u_{0,\varepsilon},$$
we compute
\[
(-\Delta+1)\int_0^\infty
e^{s(\Delta-1)}u_{0,\varepsilon}\dd s=-\int_0^\infty
\frac d{ds}\left(e^{s(\Delta-1)}u_{0,\varepsilon}\right)\dd s=u_{0,\varepsilon}.
\]
By uniqueness for \eqref{eq:concentrated-signal}, we obtain
\[
v_{0,\varepsilon}
=\int_0^\infty e^{s(\Delta-1)}u_{0,\varepsilon}\dd s.
\]
Applying the smoothing Neumann semigroup estimates to this representation,
we find, for every $q\in(1,2)$, a constant $C_{59}>0$ such that
\begin{equation}\label{v0eps-bdd}
\|v_{0,\varepsilon}\|_{1,q}
\le C_{59}\int_0^\infty \left(1+s^{-\frac32+\frac1q}\right)e^{-s}\|u_{0,\varepsilon}\|_1ds\le \frac{4C_{59}}{2-q}m.
\end{equation}
For brevity, set
$a_\varepsilon:=m-\varepsilon^2/\chi$,
and
$b_\varepsilon:=\varepsilon^2/(\chi|\Om|)$.
Then
$u_{0,\varepsilon}(x)
=a_\varepsilon\varepsilon^{-2}\varphi(x/\varepsilon)
+b_\varepsilon$
and $a_\varepsilon+|\Om|b_\varepsilon=m$.
Changing variables $x=\varepsilon y$ and using
$\int_{B_1}\varphi=1$, we obtain  $\int_\Om u_{0,\varepsilon}=m$, and
\begin{equation}\label{eq:concentrated-entropy-identity}
\int_\Om
a_\varepsilon\varepsilon^{-2}
\varphi\left(\frac{x}{\varepsilon}\right)
\ln\left[
a_\varepsilon\varepsilon^{-2}
\varphi\left(\frac{x}{\varepsilon}\right)
\right]=a_\varepsilon\left[
\ln a_\varepsilon+2\ln\frac1\varepsilon
+\int_{B_1}\varphi\ln\varphi
\right].
\end{equation}
By the concavity of $\ln (1+z)$ for $z\ge0$, we  apply the entropy inequality
\[
0
\le(X+Y)\ln(X+Y)-X\ln X-Y\ln Y\le (X+Y)\ln2,
\quad X, Y\ge0,
\]
with the convention $0\ln0=0$, to
$X=a_\varepsilon\varepsilon^{-2}\varphi(x/\varepsilon)$
and $Y=b_\varepsilon$.
Integrating the resulting inequalities over $\Om$ and using
\eqref{eq:concentrated-entropy-identity}
together with $\int_\Om u_{0,\varepsilon}=m$,
we obtain
\begin{equation}\label{eq:initial-entropy-splitting-bound}
0\le\int_\Om u_{0,\varepsilon}\ln u_{0,\varepsilon}-a_\varepsilon\left(
\ln a_\varepsilon+2\ln\frac1\varepsilon
+\int_{B_1}\varphi\ln\varphi\right)
-|\Om|b_\varepsilon\ln b_\varepsilon
\le m\ln2.
\end{equation}
By \eqref{eq:initial-entropy-splitting-bound},
the initial entropy equals the sum of the entropies
of the concentrated part and the constant background,
plus an error between $0$ and $m\ln2$.
Using $a_\varepsilon=m-|\Om|b_\varepsilon$,
we also have
\begin{equation}\label{eq:initial-entropy-remainder}
\begin{aligned}
&a_\varepsilon\left(
\ln a_\varepsilon+2\ln\frac1\varepsilon
+\int_{B_1}\varphi\ln\varphi
\right)
+|\Om|b_\varepsilon\ln b_\varepsilon
-2m\ln\frac1\varepsilon\\
={}&a_\varepsilon\ln a_\varepsilon
+a_\varepsilon\int_{B_1}\varphi\ln\varphi
-2|\Om|b_\varepsilon\ln\frac1\varepsilon
+|\Om|b_\varepsilon\ln b_\varepsilon.
\end{aligned}
\end{equation}
The right-hand side of \eqref{eq:initial-entropy-remainder}
remains bounded as $\varepsilon\to0$,
since $a_\varepsilon\to m>0$,
$b_\varepsilon\ln(1/\varepsilon)\to0$,
$b_\varepsilon\ln b_\varepsilon\to0$, and
$\int_{B_1}|\varphi\ln\varphi|<\infty$.
Combining \eqref{eq:initial-entropy-splitting-bound}
with \eqref{eq:initial-entropy-remainder},
we obtain a constant $C_{60}>0$ such that
\[
\left|
\int_\Om u_{0,\varepsilon}\ln u_{0,\varepsilon}
-2m\ln\frac1\varepsilon
\right|\le C_{60}.
\]
Consequently,
\begin{equation}\label{eq:initial-entropy}
\int_\Om u_{0,\varepsilon}\ln u_{0,\varepsilon}
=2m\ln\frac1\varepsilon+O(1),
\qquad\text{as }\varepsilon\to0+.
\end{equation}
By the well-known uniform expansion of the Green function
$G(x,y)=-(2\pi)^{-1}\ln|x-y|+O(1)$
for $x,y\in B_{R/2}$
\cite[Eq.~(1.5)]{bib:wei-winter-2002-strong-coupling}
and the identity $\int_\Om G(x,y)\dd y=1$, associated with \eqref{eq:concentrated-signal},
the Green remainder and background terms contribute $O(1)$.
Recalling the definition of $u_{0,\varepsilon}(x)=a_\varepsilon \varepsilon^{-2}\varphi(x/\varepsilon)+b_\varepsilon$ with $a_\varepsilon=m-\varepsilon^2/\chi$ and $b_\varepsilon=\varepsilon^2/(\chi|\Omega|)$
and using $\int_{B_1}\varphi=1$, we use the change of variables $x=\varepsilon \sigma$ and $y=\varepsilon \zeta$    to deduce 
\begin{equation}\label{eq:initial-interaction}
\begin{aligned}
&\chi\int_\Om u_{0,\varepsilon}v_{0,\varepsilon}\\
=&\chi\int_\Om u_{0,\varepsilon}(x)
\int_\Om G(x,y)u_{0,\varepsilon}(y)\dd y\dd x\\
=&\chi\int_\Om\int_\Om
G(x,y)\left[a_\varepsilon \varepsilon^{-2}\varphi(\frac{x}{\varepsilon})+b_\varepsilon\right]\left[a_\varepsilon \varepsilon^{-2}\varphi(\frac{y}{\varepsilon})+b_\varepsilon\right]
\dd x\dd y\\
=&\chi a_\varepsilon^2\varepsilon^{-4}\int_\Om\int_\Om
G(x,y)  \varphi(\frac{x}{\varepsilon})\varphi(\frac{y}{\varepsilon})\dd x\dd y+O(\varepsilon^2)\\
=&\chi\left(m-\frac{\varepsilon^2}{\chi}\right)^2
\int_{B_1}\int_{B_1}
G(\varepsilon\sigma,\varepsilon\zeta)
\varphi(\sigma)\varphi(\zeta)\dd\sigma\dd\zeta+O(\varepsilon^2)\\
=&\frac{\chi\left(m-\frac{\varepsilon^2}{\chi}\right)^2}{2\pi}
\int_{B_1}\int_{B_1}
\varphi(\sigma)\varphi(\zeta)
\ln\frac{1}{\varepsilon|\sigma-\zeta|}
\dd\sigma\dd\zeta+O(1)\\
=&\frac{\chi^2m^2-2\chi m\varepsilon^2+\varepsilon^4}{2\pi\chi}
\left(\ln\frac1\varepsilon+\int_{B_1}\int_{B_1}
\varphi(\sigma)\varphi(\zeta)\ln\frac1{|\sigma-\zeta|}\dd\sigma\dd\zeta\right)+O(1)\\
=&\frac{\chi m^2}{2\pi}\ln\frac1\varepsilon+O(1),
\end{aligned}
\end{equation}
in the last equality, we used the  integrability of the logarithm  function:
\begin{align*}
&\int_{B_1}\int_{B_1}
\ln\frac1{|\sigma-\zeta|}\dd\sigma\dd\zeta\\
\le& \int_{B_1}\left(\int_{B(\zeta, e^{-1})}\ln\frac1{|\sigma-\zeta|}\dd\sigma+\int_{B_1\setminus B(\zeta, e^{-1})}\ln\frac1{|\sigma-\zeta|}\dd\sigma\right)\dd\zeta\\
\le& \left(1+\frac12 e^{-2}\right)\pi^2.
\end{align*}
Testing \eqref{eq:concentrated-signal} by
$v_{0,\varepsilon}$ and integrating by parts, we derive
\begin{equation}\label{eq:initial-signal-testing}
\int_\Om
\left(|\nabla v_{0,\varepsilon}|^2+v_{0,\varepsilon}^2\right)
=\int_\Om u_{0,\varepsilon}v_{0,\varepsilon}.
\end{equation}
Using the energy definition of $\mathcal{F}$ in \eqref{eq:full-energy} and
\eqref{eq:initial-signal-testing}, and then substituting
\eqref{eq:initial-entropy} and \eqref{eq:initial-interaction},
we obtain
\[
\begin{aligned}
\mathcal F(u_{0,\varepsilon},v_{0,\varepsilon},w_0)
+\xi\int_\Om u_{0,\varepsilon}w_0&=\int_\Om u_{0,\varepsilon}\ln u_{0,\varepsilon}
-\frac\chi2\int_\Om u_{0,\varepsilon}v_{0,\varepsilon}\\
&=\frac{m(8\pi-m\chi)}{4\pi}
\ln\frac 1\varepsilon+O(1).
\end{aligned}
\]
Since $ m\chi>8\pi$, the right-hand side tends to $-\infty$ as $\varepsilon\to0$, and we obtain \eqref{eq:energy-expansion}.
\end{proof}

\begin{proof}[Proof of Theorem~\ref{thm:mass}]
We fix $B,T_0>0$ and $m>8\pi/\chi$ and choose the initial data
as in Lemma~\ref{lem:initial-energy}.
With $q=3/2$ and
$\mu=\max\{1,8C_{59}\}m\chi$, we have from \eqref{v0eps-bdd} that
$\chi\max\left\{\|u_{0,\varepsilon}\|_1,
\|v_{0,\varepsilon}\|_{1,q}\right\}
\le\mu$.
Hence the constants $L$ and $C_{57}$ in
Lemma~\ref{prop:criterion} can be chosen independently
of $\varepsilon$ and $w_0\in\mathcal W_B^0$.
For sufficiently small $\varepsilon>0$,
the asymptotic estimate \eqref{eq:energy-expansion} yields 
\[
\begin{split}
-\mathcal F(u_{0,\varepsilon},v_{0,\varepsilon},w_0)
-\xi\int_\Om u_{0,\varepsilon}w_0
=&\frac{m(m\chi-8\pi)}{4\pi}
\ln\frac{1}{\varepsilon}+O(1)\\
\ge& \frac{m(m\chi-8\pi)}{8\pi}
\ln\frac{1}{\varepsilon}\to +\infty.
\end{split}
\]
Thus \eqref{eq:thresholds} holds, and Lemma~\ref{prop:criterion} gives
\[
T_{m,\varepsilon}
\le12L\left[\frac{8\pi}{m(m\chi-8\pi)}\right]^{1/3}
\left(\ln\frac 1\varepsilon\right)^{-1/3}.
\]
Since the right-hand side tends to zero as $\varepsilon\to0$,
we can choose $\varepsilon_0>0$ sufficiently small
so that \eqref{eq:lifespan-eps} holds for all
$\varepsilon\in(0,\varepsilon_0)$ and every $w_0\in\mathcal W_B^0$.
Fixing one such $\varepsilon$ and setting
$(u_0,v_0)=(u_{0,\varepsilon},v_{0,\varepsilon})$, we obtain
\eqref{eq:blowup} for every $w_0\in\mathcal W_B^0$.
\end{proof}

\section*{Acknowledgments}
\addcontentsline{toc}{section}{Acknowledgments}
This work was partially supported by NSF of China (Nos. 12371214 and  12071476).

\end{document}